\documentclass[12pt,a4paper]{amsart}
\usepackage{amscd,mathrsfs,amssymb}
\usepackage[all]{xy}

\newtheorem{theorem}{Theorem}[section]
\newtheorem{prop}[theorem]{Proposition}
\newtheorem{corol}[theorem]{Corollary}
\newtheorem{lemma}[theorem]{Lemma}

\theoremstyle{definition}
\newtheorem{defin}[theorem]{Definition}

\theoremstyle{remark}

\numberwithin{equation}{section}

\def\iff{if and only if }

\def\tf{torsion free}
\def\oc{one-to-one correspondence}
\def\nc#1{(^#1)}
\def\nn#1{(_#1)}
\def\rr#1{\stackrel{#1}{-}}
\def\dlim{\varinjlim}

	\def\cw{\mathsf{CW}}
\def\cf{\mathsf{CF}}		\def\dg{\dagger}
	\def\pis{\pi^s}	

\def\es{\mathop\mathsf{Es}\nolimits}
\def\hos{\mathop\mathsf{Hos}\nolimits}
\def\hot{\mathop\mathsf{Hot}\nolimits}
\def\bim{\mathop\mathsf{Bim}\nolimits}
\def\mod{\mathrel\mathrm{mod}}
\def\im{\mathop\mathrm{im}nolimits}
\def\id{\mathop\mathrm{Id}nolimits}
\def\cok{\mathop\mathrm{Coker}nolimits}
\def\diag{\mathop\mathrm{diag}nolimits}
\def\Hom{\mathop\mathrm{Hom}nolimits}
\def\Mat{\mathop\mathrm{Mat}nolimits}

\def\set#1{\left\{\,#1\,\right\}}
\def\setsuch#1#2{\left\{\,#1\mid #2\,\right\}}
\def\row#1#2{\left( #1_1 , #1_2 , \dots , #1_{#2} \right)}
\def\mtr#1{\begin{pmatrix}#1\end{pmatrix}}

	\def\bee{\bigvee}
\def\pis{\pi^S}		
\def\df{\mbox{-}}	

\def\ccc{{_\cB}\cC_{\cA}}
\def\cel{*=0{\bullet}}
\def\cels{*=0{\mbox{\scriptsize$\bullet\bullet\bullet$}}}
\def\clls{*=0{\mbox{\scriptsize$\bullet\bullet\bullet\bullet$}}}
\def\lin{\ar@{-}}	
\def\sbe{\subseteq}   \def\sb{\subset}

\def\scs{\scriptsize\bf}
\def\phan{\phantom{\bigg|}}
\def\Mod{\mbox{-}\mathsf{mod}}
 \def\phh{\phantom{C^7:}}

\def\cA{\mathscr A}  \def\cB{\mathscr B}  \def\cI{\mathscr I}
\def\cJ{\mathscr J}  \def\cS{\mathscr S}  \def\cT{\mathscr T}
\def\cV{\mathscr V}  \def\cC{\mathscr C}  \def\cD{\mathscr D}
\def\al{\alpha}    \def\be{\beta}    \def\ga{\gamma}
\def\io{\iota}     \def\La{\Lambda}  \def\Ga{\Gamma}
\def\bB{\mathbf B}  \def\bS{\mathbf S}
\def\mZ{\mathbb Z}  
\def\mN{\mathbb N}  \def\mD{\mathbb D}
\def\rH{\mathrm H}  \def\fR{\mathbf r}
\def\dE{\mathfrak E}   \def\dF{\mathfrak F}
\def\dX{\mathfrak X}   \def\Mk{\Bbbk}

\def\*{\otimes}  \def\8{\infty}   \def\+{\oplus}
\def\bop{\bigoplus}   \def\xx{\times}
\def\ti{\tilde}   \def\bup{\bigcup}
\def\xarr{\xrightarrow}   \def\mps{\mapsto}

\usepackage{hyperref}

\begin{document}

\title[Matrix problems and stable homotopy]
{Matrix problems, triangulated categories and stable homotopy types}

\author{Yuriy A. Drozd}
\address{Institute of Mathematics, National Academy of Sciences of Ukraine,
 Tereschenkivska 3, 01601 Kiev, Ukraine}
\email{drozd@imath.kiev.ua, y.a.drozd@gmail.com}
\urladdr{www.imath.kiev.ua/$\sim$drozd}
\thanks{This paper summarizes the investigations mainly made during the
  stay of the author at the Max-Plank-Insitut f\"ur Mathematik in Bonn.
 The author was also partially supported by the INTAS Grant
  06-1000017-9093.}

 \begin{abstract}
  We show how the matrix problems can be used in studying triangulated
  categories. Then we apply the general technique to the
  classification of stable homotopy types of polyhedra, find out the
  ``representation types'' of such problems and give a description of
  stable homotopy types in finite and tame cases.
   \end{abstract}

 \maketitle
 \tableofcontents


The technique of matrix problems, especially, of bimodule categories,
has proved their efficiency in lots of problems from representation
theory, algebraic geometry, group theory and other domains of modern
algebra. During last years, mainly due to the works of Baues, Henn,
Hennes, and the author, it has found new applications in algebraic
topology, namely, in studying stable homotopy classes of polyhedra
(see \cite{hen}, \cite{bh}, \cite{bd1}--\cite{bd4}, \cite{d2}). In the
survey \cite{d1} the author has picked out the background of this
approach, which is based on the trianguled structure of thew stable
homotopy category. In this paper we show that the same method can be
used in general situation, when we construct subcategories of a
triangulated category from simpler ones (see Section~1). Then we
summarize what can be done using this technique for the classification
problem of stable homotopy classes. Namely, we consider the
subcategories $\cS_n$ of the stable homotopy category consisting of
polyhedra having only cells of $n$ consecutive dimensions. We classify
polyhedra from $\cS_n$ for $n\le 4$ and show that for $n>4$ this
problem becomes \emph{wild} in the sense of the representation theory
of algebras. Then we consider the subcategories $\cT_n$ of $\cS_n$
consisting of polyhedra with no torsion in ontegral homologies. This
time we classify polyhedra from $\cT_n$ for $n\le7$ and show that for
$n>7$ their classification is also a wild problem. In some sense,
these results are ``final,'' though we are sure that this technique
will be useful for some other problems of algebraic topology as well
as for studying other triangulated categories.

Since the technical details of calculations are sometimes rather
cumbersome and can be found in the previous papers, we usually omit
them, just outlining the ideas.

 \section{Matrix problems arising in triangulated categories}
 \label{s1}

 Let $\cC$ be a triangulated category with the shift $A\mps SA$, $\cA$ and $\cB$
 be two fully additive (but usually not triangulated) subcategories of $\cC$. We denote by
 $\cA\dg\cB$ the full subcategory of $\cC$ consisting of all objects $C$ arising in triangles 
 \begin{equation}\label{e11}
  A \xarr a B \xarr b C \xarr c SA \ \text{ with }\ A\in\cA,\ B\in\cB.
 \end{equation}
 We also denote by $\cI$ the ideal of the category $\cC$ consisting of all morphisms
 $\ga:C\to C'$ that factorizes both through $\cB$ and through $S\cA$, i.e. such that
 $\ga=\ga'\al=\ga''\be$, where $\al:C\to SA,\,\be:C\to B$, where $A\in\cA,\, B\in\cB$.

 On the other hand, we consider the $\cA\df\cB$-bimodule $\ccc$, which is
 the restriction of the regular $\cC$-bimodule $\cC(A,B)$ for $A\in\cA,\,B\in\cB$.
 We often omit subscripts and denote this bimodule by $\cC$ if it cannot lead to
 misunderstanding. Recall that the \emph{bimodule category} $\bim(\ccc)$
 has $\bup_{\substack{A\in\cA\\B\in\cB}} \cC(A,B)$ as the set of objects, while
 the set of morphisms $\bim(a,a')$, where $a:A\to B,\,a':A'\to B'$, is defined as 
 \[
  \setsuch{(\al,\be)}{\al:A\to A',\,\be:B\to B',\ \be a=a'\al}.
 \]
 We denote by $\cJ$ the ideal of $\bim(\ccc)$ consisting of all morphisms $(\al,\be):a\to a'$
 such that $\al$ factors through $a$ and $\be$ factors through $a'$.

 We define a functor $F:\bim(\ccc)\to(\cA\dg\cB)/\cI$ as follows. For every morphism
 $a:A\to B$, choose a triangle like \eqref{e11} and set $C=Fa$. If $a':A'\to B',\ C'=Fa'$
 and $(\al,\be)\in\bim(a,a')$, there is $\ga:C\to C'$ such that the diagram
 \begin{equation}\label{e12}
  \begin{CD}
   A @>a>> B  @>b>> C @>c>> SA @>Sa>> SB \\
   @V\al VV	 @V\be VV   @V\ga VV  @VVS\al V  @VVS\be V \\
   A' @>a'>> B'  @>b'>> C' @>c'>> SA' @>Sa'>> SB' 
 \end{CD}
 \end{equation} 
 commutes. Set $F(\al,\be)=\ga\mod\cI$. We must check that the latter definition is
 consistent. Indeed, if $\ga':C\to C'$ is another morphism making diagram \eqref{e12}
 commutative, $g=\ga-\ga'$, then $gb=c'g=0$, therefore there are $f:SA\to C'$ and
 $h:C\to B'$ such that $g=cf=b'h$, i.e. $g\in\cI$. Thus $F$ is well-defined. 

 Suppose now that $\cC(B,SA)=0$ for all $A\in\cA,\,B\in \cB$. In this situation
 we define a functor $G:\cA\dg\cB\to\bim(\ccc)/\cJ$ as follows. Let $C\in\cA\dg\cB$.
 Choose one triangle like \eqref{e11} and set $a=GC$. If $GC'=a'$, i.e.
 $C'$ occur in the triangle
 \[
     A' \xarr {a'} B' \xarr {b'} C' \xarr {c'} SA' \ \text{ with }\ A'\in\cA,\ B'\in\cB,
 \]
 and $\ga:C\to C'$, then $c'\ga b=0$, hence $\ga b=b'\be$ for some $\be:B\to B'$.
 Choose one of such triangles
 Since 
 \[
  B \xarr b C \xarr c SA \xarr{-Sa} SB \\
 \]
 and
 \[
  B' \xarr b' C' \xarr c' SA' \xarr{-Sa'} SB' 
 \]
 are also triangles, there is a morphism $\al:A\to A'$ that makes the diagram \eqref{e12}
 commutative, thus $(\al,\be)\in\bim(a,a')$. Set $G\ga=(\al,\be)\mod\cJ$. If $(\al',\be')$
 is another pair making \eqref{e12} commutative, then $(\be-\be')b'=0$, hence $\be-\be'=a'f$
 for some $f:B\to A'$; in the same way $S\al-S\al'=g(Sa)$, i.e. $\al-\al'=g'a$ for some
 $g:SB\to SA'$ and $g':B\to A'$ such that $Sg'=g$. Therefore $(\al-\al',\be-\be')\in\cJ$,
 so the functor $G$ is well-defined.
 
 \begin{theorem}\label{11}
 Suppose that $\cC(B,SA)=0$ for all $A\in\cA,\,B\in\cB$. Then
  the functors $F,G$ constructed above induce quasi-inverse functors
 $\bar F:\bim(\ccc)/\cJ\to(\cA\dg\cB)/\cI$ and $\bar G:(\cA\dg\cB)/\cI\to\bim(\ccc)/\cJ$.
 Thus $(\cA\dg\cB)/\cI\to\bim(\ccc)/\cJ$. Moreover, $\cI^2=0$, therefore, the natural
 functor $\Pi:(\cA\dg\cB)\to(\cA\dg\cB)/\cI$ is an epivalence.
 \end{theorem}

 Recall that an \emph{epivalence} is a functor $E:\cC_1\to\cC_2$, which is
 \begin{itemize}
\item  \emph{full}, i.e. all induced maps $\cC_1(X,Y)\to\cC_2(EX,EY)$ are surjective;
 \item  \emph{dense}, i.e. every object from $\cC_2$ is isomorphic to $EX$ for some $X\in\cC_1$;
 \item  \emph{conservative}, i.e. $f\in\cC_1(X,Y)$ is invertible \iff so is $Ef\in\cC_2(EX,EY)$.
\end{itemize}
 (In \cite{ba1} such functors are called \emph{detecting}.)
 Note that then also
 \begin{itemize}
  \item  $X\simeq Y$ in $\cC_1$ \iff $EX\simeq EY$ in $\cC_2$;
\item  if $\cC_1$ and $\cC_2$ are additive, then an object $X\in\cC_1$ is indecomposable (into a
 nontrivial direct sum) \iff so is $EX$.
\end{itemize}
 
 \begin{proof}
  One immediately sees that $F(\cJ)=0$ and $G(\cI)=0$, hence $\bar F$ and $\bar G$ are well-defined.
 Moreover, we have already seen that, given $(\al,\be)$, the morphism $\ga$ is defined up to a summand
 from $\cI$, and given $\ga$, the pair $(\al,\be)$ is defined up to a summand from $\cJ$. It
 obviously implies that $\bar F\bar G\simeq\id$ and $\bar G\bar F\simeq\id$. If $\ga:C\to C'$ and
 $\ga':C'\to C''$ are from $\cI$, then $\ga=gf$ for some $f:C\to  B$ and $g:B\to C'$, where
 $B\in\cB$, while $\ga'=g'f'$ for some $f':C'\to SA$ and $g':A\to C''$, where $A\in\cA$.
 Then $\ga'\ga=g'f'gf=0$, since $f'g\in\cC(B,SA)$. Thus $\cI^2=0$ and, therefore,
 $\Pi$ is an epivalence.
 \end{proof} 

  \begin{corol}\label{12}
  Under conditions of Theorem~\ref{11}, let $\cV$ be a subbimodule of $\ccc$ such that
 $f_1af_2=0$ whenever $a\in\cV,\,f_i\in\cC(B_i,A_i)$ with $A_i\in\cA,\,B_i\in\cB$
 $(i=1,2)$. Denote by $\cA\dg_\cV\cB$ the full subcategory of $\cA\dg\cB$ consisting
 of all objects $C$ arising in triangles \eqref{e11} with $a\in\cV$,
 $\cI_\cV=\cI\cap(\cA\dg_\cV\cB)$, $\cJ_\cV=\cJ\cap\bim(\cV)$. Then the functor
 $F$ and $G$ constructed above induce quasi-inverse functors 
 $\bar F:\bim(\cV)/\cJ\to(\cA\dg_\cV\cB)/\cI_\cV$ and
 $\bar G:(\cA\dg_\cV\cB)/\cI\to\bim(\cV)/\cJ_\cV$.
 Thus $(\cA\dg_\cV\cB)/\cI\simeq\bim(\cV)/\cJ_\cV$. Moreover, $\cI_\cV^2=0$
 and $\cJ_\cV^2=0$, therefore, the natural functors $(\cA\dg_\cV\cB)\to(\cA\dg_\cV\cB)/\cI_\cV$
 and $\bim(\cV)\to\bim(\cV)/\cJ^2_\cV$ are epivalences. In particular, there is a \oc s 
 between isomorphism classes of objects and of indecomposable objects from $\cA\dg_\cV\cB$ and $\bim(\cV)$.
 \end{corol}

 \section{Stable homotopy category}
 \label{s2}

 In this paper the word ``\emph{polyhedron}'' is used as a synonym for``\emph{finite cell}
 (or CW) \emph{complex}''. We denote by $\hot$ the category of punctured topological spaces
 with homotopy classes of continuous maps as morphisms and by $\cw$ its full subcategory
 consisting of polyhedra. We denote by $CX$ the \emph{cone over the space} $X$, i.e. the
 factor space $X\xx I/X\xx 1$, $I=[0,1]$ being the unit interval. For a map $f:X\to Y$ we
 denote by $Cf$ the \emph{cone of this map}, i.e. the factor space $(Y\sqcup CX)/\sim$,
 where the equivalence relation $\sim$ is given by the rule $f(x)\sim (x,0)$. Let also $SX$
 be the \emph{suspension} of $X$, i.e. the factor space $CX/(X\xx 0)$. This operation
 induces a functor $S:\hot\to\hot$. Note that for every $X$ the space $SX$ is an
 \emph{$H$-cogroup} and the $n$-fold suspension $S^nX$ is a \emph{commutative
 $H$-cogroup} for $n\ge2$ \cite[2.21\,--\,2.26]{sw}. Therefore, $\hot(S^nX,Y)$
 is a group, commutative for $n\ge2$. The natural maps
 $\hot(S^nX,S^nY)\to\hot(S^{n+1}X,S^{n+1}Y)$ are group homomorphisms. Set
 \[
  \hos(X,Y)=\dlim_n \hot(S^nX,S^nY).
 \]
 It is a group called the group of \emph{stable maps} from $X$ to $Y$. Thus we get the
 \emph{stable homotopy category} $\hos$ and its full subcategory $\cS$ consisting of
 polyhedra. We also denote by $\cf$ and $\cT$ respectively the full subcategories of $\cw$
 and of $\cS$ consisting of \emph{torsion free polyhedra} $X$, i.e. such that all integral
 homology groups $\rH_k(X)=\rH_k(X,\mZ)$ are torsion free. The groups $\hos(S^n,X)$
 are called the \emph{stable homotopy groups} of the space $X$ and denoted by $\pis_n(X)$.
 
 The category $\hos$ is additive, with the
 \emph{bouquet} (or \emph{wedge}) $X\vee Y$ playing the role of direct sum. Moreover,
 $\hos$ is \emph{fully additive}, i.e. every idempotent in it splits \cite[Theorem~4.8]{co}.
 The suspension induces a functor, which we also denote by $S:\hos\to\hos$. Obviously,
 it is fully faithful. Thus we can ``supplement'' it so that $S$ becomes an equivalence. To do it,
 we consider formal ``imaginary spaces'' $S^nX$ with $n<0$ setting, for $n<0$ or $m<0$,
 $\hos(S^nX,S^mY)= \hos(S^{n+k}X,S^{m+k}Y)$, where $k=-\min(n,m)$. Then we consider
 formal bouquets $\bee_{i=1}^rX_i$, where each $X_i$ is either a ``real'' or an ``imaginary''
 space, and define $\hos(\bee_{j=1}^sY_j,\bee_{i=1}^rX_i)$ as the set of $r\xx s$
 matrices $(f_{ij})$ with $f_{ij}\in\hos(Y_j,X_i)$ (see \cite{co} for details). As a result we
 get the category (also denoted by $\hos$), where $S$ is an auto-equivalence.

 In fact, the new category is a \emph{triangulated} category. The triangles in it are the
 \emph{cofibration sequences}, i.e. those isomorphic to the \emph{cone sequences}
 \[
  X\xarr f Y \xarr g Cf \xarr h SX,
 \]
 where $g$ is the natural embedding $Y\to Cf$ and $h$ is the natural surjection
 $Cf\to SX\simeq Cf/Y$ \cite{pu}. Note that in the stable category $\hos$ they
 coincide with the \emph{fibration sequences} \cite{co}, though we do not use this fact.

 We denote by $\cw^k_n$ the full subcategory of $\cw$ consisting of $(n-1)$-connected
 cell complexes of dimension at most $n+k$. If $X\in\cw^k_n$, one can suppose that
 its $(n-1)$-th skeleton $X^{n-1}$ (the ``$(n-1)$-dimensional part'' of $X$) consist of a
 unique point and $X$ has no cells of dimensions greater than $n+k$. Following Baues,
 we describe such a cell complex using its \emph{gluing} (or \emph{attachment})
 \emph{diagram}, which looks like (for $n=7,k=6$)
 \begin{equation}\label{e21}
     \vcenter
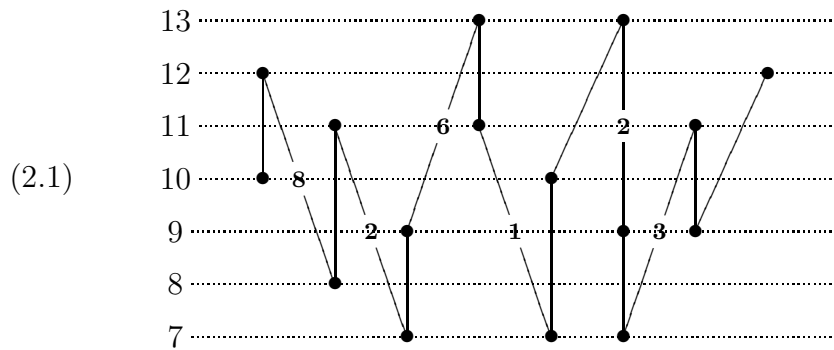
{\xymatrix@R=1ex@C=1.7em{
   13 \ar@{.}[rrrrrrrrr] &&&& \cel \ar@{-}[dd]
   \ar@{-}[ddddl]|*+\txt{\scs 6} && \cel \ar@{-}[dddd]|*+\txt{\scs2} 
     \ar@{-}[dddl] &&& *=0{}\\
   12 \ar@{.}[rrrrrrrrr] & \cel \ar@{-}[dd]
   \ar@{-}[ddddr]|*+\txt{\scs8} &&&&&&& \cel
   \ar@{-}[dddl] & *=0{}\\ 
   11 \ar@{.}[rrrrrrrrr] && \cel \ar@{-}[ddd]
   \ar@{-}[ddddr]|*+\txt{\scs2} && \cel \ar@{-}[ddddr]|*+\txt{\scs1}
   &&& \cel 
   \ar@{-}[ddddl]|*+\txt{\scs3} \ar@{-}[dd] && *=0{}\\
   10 \ar@{.}[rrrrrrrrr] & \cel &&&& \cel \ar@{-}[ddd] &&&& *=0{}\\
   9 \ar@{.}[rrrrrrrrr]  &&&\cel \ar@{-}[dd] &&& \cel \ar@{-}[dd] & \cel && *=0{}\\
   8 \ar@{.}[rrrrrrrrr] && \cel &&&&&&& *=0{}\\
   7 \ar@{.}[rrrrrrrrr] &&& \cel && \cel & \cel  &&& *=0{}	}} 
 \end{equation}%
   In this diagram each bullet on the level $m$ corresponds to an $m$-dimensional cell, i.e. to a ball
  $\bB^m$ glued to the $(m-1)$-dimensional skeleton $X^{m-1}$  by a map of its boundary
 $f:\bS^{m-1}\to X^{m-1}$. The lines between this bullet  and the lower ones describe the
 nonzero components of the map $f$. If there are more than
  one nonzero map between $\bS^{m-1}$ and a smaller $\bS^l\ (l<m)$,
  these lines carry some marks precising the corresponding maps. Especially, in our example the groups
 $\hos(S^{l+3},S^l)$ are cyclic of order $24$, so we put the marks that show, which multiple
 of the generator is used for this gluing. There are no marks on other lines, since the groups
 $\hos(S^{l+2},S^l)$ are of order $2$, so only have one nonzero element. 

 Every polyhedron from $\cS$ decomposes into a direct sum of indecomposable ones. Note that
 such a decomposition is far from being unique (see \cite[4.2]{co} for examples). Nevertheless, a
 description of indecomposable polyhedra in $\cS$ can be a good first step towards the
 classification problem. Moreover, if the endomorphism ring $\es(Y)=\hos(Y,Y)$ is \emph{local}
 and $Y\vee Y'\simeq\bop_iY_i$, there is an index $i$ such that $Y_i\simeq Y\vee Y''$
 \cite[Lemma~I.3.5]{bas} 
 Hence, in all decompositions of a polyhedron $X$ into bouquets of indecomposables the multiplicity
 of $Y$ is the same. Another approach gives the notion of
 \emph{congruence}. Namely, we say that two polyhedra $X,Y$ are
 \emph{congruent} if there is a polyhedron $Z$ such that $X\vee
 Z\simeq Y\vee Z$. One can show, following \cite{jac} or \cite{gen},
 that an equivalent condition is that 
 the images of $X$ and $Y$ in all \emph{localizations} $\cS_p$ of the
 stable homotopy category are isomorphic. Here $\cS_p$ ($p$ is a
 prime integer) is the category whose objects are polyhedra, but
 $\hos_p(X,Y)=\hos(X,Y)\*\mZ_p$, where $\mZ_p$ is the ring of
 $p$-adique integers. (The same notion is obtained if we replace
 $\mZ_p$ by the subring $\setsuch{a/b}{a,b\in\mZ,\,p\nmid b}$ of the
 rational numbers.) We call the classes of congruence \emph{genera},
 like they do in the theory of integral representations. Though genera
 satisfy the cancellation property (in fact, by definition), their
 decomposition into bouquets of indecomposable is not unique too (see
 the first of the cited examples from \cite{co}).

 Recall that due to the Generalized Freudenthal Theorem \cite[Theorem~1.21]{co} there is no need to
 go up to infinity in defining $\hos(X,Y)$ if we deal with polyhedra. Namely, if $Y$  is $(n-1)$-connected
 and $\dim X\le m$, then the map $\hot(X,Y)\to\hot(SX,SY)$ is bijective if $m<2n-1$ and
 surjective if $m=2n-1$. It implies that the map $\hot(S^kX,S^kY)\to\hos(X,Y)$ is bijective for
 $k>m-2n+1$ and surjective for $k=m-2n+1$. In particular, if $Y$ is $(n-1)$-connected,
 $\pis_m(Y)\simeq\pi_{2(m-n+1)}(S^{m-n+2}Y)$.
 Moreover, on the subcategory of $\hot$ consisting of simply connected
 spaces the suspension functor is conservative. Therefore, the induced functor $\cw^k_n\to\cw^k_{n+1}$
 is an equivalence for $n>k+1$ and an epivalence for $n=k+1$. Denote by $\cS_n$ the image in
 $\cS$ of the category $\cw^{n-1}_n$. The polyhedra from $\cS_n$ can only have cells on
 $n$ consecutive levels (from $n$-th up to $(2n-1)$-th) and every polyhedra having cells on $n$
 consecutive levels is isomorphic in $\cS$ to $S^mX$ for some integer $m$ and some $X\in\cS_n$.
 We also denote by $\cT_n$ the full subcategory of $\cS_n$ consisting of torsion free polyhedra.
 
  \begin{defin}[cf.~\cite{ba}]\label{21}
  An \emph{atom} is an indecomposable object $A$ from $\cS_n$, which does not belong to
 $S(\cS_{n-1})\cup S^2(\cS_{n-1})$. (In other words, any polyhedron isomorphic
 to $A$ in $\cS$ must have cells of dimensions $n$ and $2n-1$.) If $A$ is an atom, all
 polyhedra of the sort $S^mA$ are called \emph{suspended atoms}.
 \end{defin}

 This definition immediately implies that every polyhedron is isomorphic in $\cS$ to a bouquet of
 suspended atoms, though, as we have mentioned, such a decomposition is not unique. 
 Note that, unlike Baues, we consider $S^1$ as an atom (a unique atom in $\cS_1$), hence all
 spheres are  considered as suspended atoms. Note also that this definition implies that all atoms 
 are of odd dimensions: an atom from $\cS_n$ is of dimension $2n-1$.

 To clarify the structure of $\cS_n$ we use the technique from Section~1. Namely, choose an
 integer $m$ such that $0\le m<n-1$ and set
  \begin{equation}\label{e22}
  \begin{split}
   \cA&=\cA_{n,m}=S^{2m+1}\cS_{n-m-1},\\ \cB&=\cB_{n,m}=S^{n-m-1}\cS_{m+1},\\ 
  \cS_{n,m}&=_{\cB_{n,m}}\!\cS\!_{\cA_{n,m}},\\
  \cI_{n,m}&=\{\,f:X\to Y\,|\,X,Y\in\cS_n,\ f \text{ factors both}\\
   &\hskip7.5em \text{through $\cB$ and  through } S\cA\},\\
  \cJ_{n,m}&=\{\,(\al,\be)\in\bim(a,a')\,|\,a,a'\in\cS_{n,m},\\ &\qquad\al \text{ factors through } a
	 \text{ and } \be \text{ factors through } a'\}.
  \end{split}
  \end{equation}
 Then polyhedra from $\cA$ only have cells in dimensions from $n+m$ up to $n-m-2$,
 while those from $\cB$ only have cell in dimensions from $n$ up to $n+m$. If $C\in\cS_n$,
 its $(n+m)$-th skeleton $B$ belongs to $\cB$, while the factor space $C/B$ belongs to $S\cA$,
 i.e. $C/B\simeq SA,\ A\in\cA$. Then $C\in\cA\dg\cB$, since $A\to B\to C\to C/B\simeq A$
 is a cofibration sequence. On the other hand, any object from $\cA\dg\cB$ obviously belongs to
 $\cS_n$. So we have proved

 \begin{theorem}\label{22}
  $\cS_n\simeq\cA_{n,m}\dg\cB_{n,m}$.
 Thus $\,\cS_n/\cI_{n,m}\simeq\bim(\cS_{n,m})/\cJ_{n,m}$.
 Moreover, $\cI_{n,m}^2=0$.
 \end{theorem}
 
 To consider torsion free polyhedra, we set
 \begin{equation}\label{e23}
   \begin{split}
     \cA^0&=\cA^0_{n,m}=S^{2m+1}\cT_{n-m-1},\\ \cB^0&=\cB^0_{n,m}=S^{n-m-1}\cT_{m+1},\\ 
  \cS^0&=\cS^0_{n,m}=\setsuch{a\in\cS_{n,m}}{\rH_{n+m}(a)=0},\\
  \cI^0_{n,m}&=\cI_{n,m}\cap(\cA^0\dg_{\cS^0}\cB^0),\\
  \cJ^0_{n,m}&=\cJ_{n,m}\cap\bim(\cS^0).
 \end{split}
 \end{equation}

 To get an analogue of Theorem~\ref{22} we need the following lemma.

 \begin{lemma}\label{23}
  Let $f\in\hos(A,B)$, where $A$ and $B$ are \tf\ polyhedra, $A$ is $(m-1)$-connected, $\dim B\le m$
 and $Cf$ is also \tf. There are decompositions $A\simeq C\+A',\ B\simeq C\+B'$ such that, with respect
 to this decomposition, $f=\mtr{\id&0\\0&g}$ with $\rH_m(g)=0$ and $Cf\simeq Cg$.
 \end{lemma}
 \begin{proof}
 Note first that if $A=kS^m,\ B=lS^m$ are bouquets of $m$-dimensional spheres, then $\rH_m(A)=m\mZ,\ 
 \rH_m(B)=l\mZ$, and the natural map $\hos(A,B)\to\Hom(\rH_m(A),\rH_m(B))$ is an isomorphism.
 In particular, every decomposition of $\rH_m(A)$ arises from a decomposition of $A$, and the same
 is true for $B$. In this case $\rH_m(f)$, or, the same, $f$ is actually an integer matrix and
 there are decompositions $A\simeq C\+A',\ B\simeq C\+B'$ (all summands are, of course,
 also bouquets of spheres) such that, with respect to them, $f=\mtr{\id&0\\0&d}$,
 where $d:A'\to B'$ can be presented by diagonal matrix without unit components.

 In general case, the calculation of homologies of cell spaces from \cite[Chapter 10]{sw} shows that the
 embedding $\al:A^m\to A$ induces a surjection $\rH_m(A^m)\to\rH_m(A)$, while the surjection
 $\be:B\to \ti B=B/B^{m-1}$ induces an embedding $\rH_m(B)\to\rH_m(\ti B)$ with torsion free
 cokernel. Therefore, there are decompositions $A^m\simeq A_1\+A_0,\ \ti B\simeq B_1\+B_0$
 such that the restriction of $\rH_m(\al)$ onto $A_1$ is an isomorphism, and that onto $A_0$ is $0$,
 while $\rH_m(f)$ induces an isomorphism $\rH_m(B)\to\im\rH_m(\be)=\rH_m(B_1)$. Denote by
 $\al_1:A_1\to A$ and $\be_1:B\to B_1$ the corresponding components of $\al$ and $\be$.
 As above, there are decompositions $A_1\simeq C\+A_0,\ B_1\simeq C\+B_0$ such that,
 with respect to them, the morphism $\be_1f\al_1=\mtr{\id&0\\0&d}$, where
 $d$ can be presented by diagonal matrix without unit components. Denote by $\io:C\to A_1$
 the natural embedding (presented by the matrix $\mtr{\id\\0}$) and by $\pi:B_1\to C$ the
 natural projection (presented by the matrix $\mtr{\id&0}$\,). Then $\pi\be_1f\al_1\io=\id$,
 so $B\simeq C\+B',\ A\simeq C\+A'$, so that, with respect to these decompositions,
 $f=\mtr{\id&0\\0&g}$. Then $Cf\simeq Cg$ and $d=\be_0g\al_0$, where
 $\al_0:A_0\to A'$ and $\be_0:B'\to B_0$. Note that $\al_0$ and $\be_0$ also induce
 isomorphisms of the $m$-th homology groups, so $\cok\rH_m(g)\simeq\cok\rH_m(d)$.
 Since this cokernel embeds in $\rH_m(Cg)$, it is torsion free. Therefore, $d=0$, whence
 $\rH_m(g)=0$.
 \end{proof}

 \begin{theorem}\label{24}
    $\cT_n\simeq\cA^0_{n,m}\dg_{\cS^0}\cB^0_{n,m}$.
 Thus $\cT_n/\cI^0_{n,m}\simeq\bim(\cS^0_{n,m})/\cJ^0_{n,m}$.
 Moreover, $(\cI^0_{n,m})^2=(\cJ^0_{n,m})^2=0$, so this equivalence induces \oc s
 between isomorphism classes of objects and of indecomposable objects in $\cT_n$
 and in $\bim(\cS^0_{n,m})$.
 \end{theorem}
 \begin{proof}
  Let $C\in\cT_n$, $B=C^{n+m}$, $SA\simeq C/B$. The triangle
 \begin{equation}\label{e24}
  A\xarr a B\xarr b C\xarr c SA
  \end{equation}
  gives rise to the exact sequence of homologies
 \begin{align*}
 &\dots \to\rH_k(A)\xarr{\rH_k(a)} \rH_k(B) \xarr{\rH_k(b)} \rH_k(C)\xarr{\rH_k(c)} \rH_k(SA)\simeq
 \\ &\simeq  \rH_{k-1}(A)\xarr{\rH_{k-1}(a)} \rH_{k-1}(B)\xarr{\rH_{k-1}(b)} \rH_{k-1}(C)\to\dots
 \end{align*}
 If $k<n+m$, then $\rH_k(A)=\rH_{k-1}(A)=0$, so $\rH_k(B)\simeq\rH_k(C)$ is torsion free.
 If $k>n+m$, we get in the same way that $\rH_k(A)\simeq\rH_{k+1}(C)$ is also torsion free. Let now
 $k=n+m$, then we get the exact sequence
  \[
 0\to \rH_{n+m+1}(C)\to \rH_{n+m}(A)\xarr{\rH_{n+m}(a)} \rH_{n+m}(B)\to\rH_{n+m}(C)\to 0.
 \]
 Note that $\rH_{n+m}(B)$ is always torsion free, since $B$
 contains no cells of dimensions bigger than $n+m$, hence $B\in\cT$. Therefore, $\rH_{n+m}(A)$
 is torsion free too, so $A\in\cT$.
 Moreover, $\cok\rH_{n+m}(a)$ is also torsion free. As both $\rH_{n+m}(A)$ and $\rH_{n+m}(B)$ are free, 
 it means that $\rH_{n+m}(A)\simeq M\+M',\ \rH_{n+m}(B)\simeq M\+M''$ so that $\rH_{n+m}(a)$   \item  
 induces isomorphism $M\to M$ and is zero on $M'$. By Lemma~\ref{23}, there are decompositions
 $A\simeq A_0\vee A',\ B\simeq A_0\vee B'$ such that, with respect to them, $a=\mtr{\id&0\\0&a'}$,
 where $a'\in\cS^0$. Then $Ca'\simeq Ca\simeq C$, so $C\in\cA^0\dg_{\cS^0}\cB^0$. On the
 other hand, if $C\in\cA^0\dg_{\cS^0}\cB^0$, i.e. belongs to a triangle \eqref{e24} with
 $A\in\cA^0,\,B\in\cB^0$ and $\rH_{n+m}(a)=0$, the exact sequence of homologies implies that
 $C\in\cT_n$. 

 To prove the remaining assertions, it is enough to show that $uav=0$ for every $a\in\cS^0(A,B)$,
 $v:B'\to A,\ u: B\to A'$, where $A,A'\in\cA^0,\ B,B'\in\cB^0$ (see Corollary~\ref{12}).
 Since $\rH_{n+m}(a)=0$, the induced map $A^{m+n}\to B/B^{m+n-1}$ is zero. On the other hand,
 $\cS(B',A/A^{m+n})=0=\cS((B')^{m+n-1},A^{m+n})$, so the map $v:B'\to A$ factors through a map
 $B'/(B')^{m+n-1}\to A^{m+n}$. Since the same holds for $u$, it implies that $uav=0$.
 \end{proof}
 
 We shall also use the following obvious lemma.

 \begin{lemma}\label{25}
  Let $X\in\cS_n$, $H_i=\rH_i(X)$. If $X$ is decomposable, there are decompositions 
 $H_i=H'_i\+H^{''}_i$ and indices $j,k$ such that both $H'_j\ne0$ and $H''_k\ne0$. 

 \emph{(Note also that $\rH_i(X)=0$ for $i<n$ or $i>2n-1$.)}
 \end{lemma}

 \section{Discrete case: Whitehead--Chang Theorem}
 \label{s3}

 We apply now Theorem~\ref{22} to polyhedra from $\cS_n$ for small $n$.
 First, we recall some values of stable homotopy groups \cite[Sections~XI.15--16]{hu}:
 \begin{itemize}
\item  $\pis_{n+1}(S^n)\simeq\mZ/2$, the generator being the (suspended) \emph{Hopf map}
 $\eta=2^{n-1}h_2$, where $h_2$ is the Hopf fibration $S^3\to S^2$;
 \item  $\pis_{n+2}(S^n)\simeq\mZ/2$, the generator being the \emph{double Hopf map}
 $\eta^2$, i.e. the composition of Hopf maps $S^{n+2}\to S^{n+1}\to S^n$;
 \item  $\pis_{n+3}(S^n)\simeq\mZ/24$, the generator being $\nu=S^{n-4}h_4$, where
 $h_4$ is the Hopf fibration $S^7\to S^4$. Moreover, the composition
 $\eta^3:S^{n+3}\to S^{n+2}\to S^{n+1}\to S^n$ equals $12\nu$.
\end{itemize}

 If $n=1$, the only atom in $\cS_1$ is $S^1$, and every polyhedron is a bouquet
 of several copies of $S^1$. If $n=2$, $\cS_2=\cA_{2,0}\dg\cB_{2,0}$, and
 $\cA_{2,0}=\cB_{2,0}=S\cS_1$. Thus every polyhedron $C$ from $\cS_2$ is isomorphic
 to the cone of a map $a:kS^2\to lS^2$.Since $\hos(S^2,S^2)=\mZ$, the map $a$
 can be considered as a matrix $(a_{ij})\in\Mat(l\xx k,\mZ)$. If $a'$ is another object
 from $\cC_{2,1}$, also considered as a matrix from $\Mat(l'\xx k',\mZ)$, a morphism
 $a\to a'$ in $\bim(\cC_{2,1})$ is given by a pair of matrices $\al\in\Mat(k'\xx k,\mZ)$,
 $\be\in\Mat(l'\xx l,\mZ)$ such that $a'\al=\be a$. Especially, this morphism is an
 isomorphism \iff both $\al$ and $\be$ are invertible. So  the well-known Smith Theorem 
 implies that every object $a\in\cC_{2,1}$ is isomorphic to one presented by a diagonal
 matrix $\diag\row qr$. Hence, every polyhedron from $\cS_2$ is isomorphic to a
 bouquet of cones $\bee_i Cq_i$, where we identify an integer $q$ with the corresponding
 map $S^2\to S^2$. Moreover, if $q=uv$, where $\gcd(u,v)=1$, then
 \[
  Cq\simeq C\mtr{1&0\\0&q}\simeq C\mtr{u&0\\0&v}\simeq Cu\vee Cv. 
 \]
 Therefore, $Cq$ can only be indecomposable if $q=p^s$, where $p$ is prime. On the other
 hand, the exact sequence of homologies arising from the triangle
 \begin{equation}\label{e31}
   S^2\xarr q S^2\xarr{} M^3(q)\xarr{} S^3 
 \end{equation}
 that $\rH_2(Cq)\simeq\mZ/q$ and $\rH_3(Cq)=0$. Hence, Lemma~\ref{25} implies that
 $Cq$ is indecomposable. Therefore, the atoms in $\cS_2$ are just $Cq$ for $q=p^s$ with
 a prime $p$. These atoms are denoted by $M^3(q)$ and their suspensions $S^kM^3(q)$ by
 $M^{k+3}(q)$. The atoms and suspended atoms $M^d(q)$ are called \emph{Moore spaces}%
 \footnote{In \cite[Section~XI.10]{hu} they are denoted by $P^d_q$ and
 called \emph{pseudo-projective spaces}.}
 \cite{co}. We also write $M^d_s$ instead of $M^d(2^s)$ (these atoms play a special role later). 

 We can calculate the groups $\hos(M^3(q),M^3(q'))$. Since $\pis_3(S^2)\simeq\mZ/2$
 \cite[Theorem~15.1]{hu}, the exact sequences for the functor $\hos$ arising from the triangles
 \eqref{e31} for $q$ and $q'$ imply that
 \begin{align} 
 &\hos(S^2,M^3(q)) \simeq \hos(M^3(q),S^3) \simeq\mZ/q,  \notag \\
 &\hos(S^3,M^3(q)) \simeq \hos(M^3(q),S^2) \simeq
 \begin{cases}
  \mZ/2 &\text{ if } q \text{ is even},\\
  0	&\text{ if } q \text{ is odd},
 \end{cases}  \notag\\
 &\hos(M^3(q),M^3(q')) \simeq \mZ/(q,q') \ \text{ if $q$ or $q'$ is odd},  \notag\\
 \intertext{and there is an exact sequence}  \label{e32}
 &\hskip2em 0 \to \mZ/2 \to \hos(M^3_s,M^3_r) \to \mZ/2^m \to 0, \text{ where } m=\min(r,s).
 \end{align}
 Note that the endomorphism rings $\es(M^3(q))$ are finite, hence, local. These considerations
 immediately imply the description of polyhedra from $\cS_2$.

 \begin{theorem}\label{31}
  Every polyhedron from $\cS_2$ uniquely (up to permutation of summands) decomposes into a 
 bouquet of spheres $S^2,S^3$ and Moore atoms $M^3(q)$.
 \end{theorem}
 
 We also need the following fact. 

 \begin{prop}\label{32}
 \[ 
 \pis_4(M^3(q))\simeq
 \begin{cases}
  0 &\text{\em if $q$ is odd},\\
  \mZ/4 &\text{\em if } q=2,\\
  \mZ/2\+\mZ/2 &\text{\em if } q=2^s,\ s>1. 
 \end{cases}
 \]
 \end{prop}  
 \begin{proof}
  Recall that $\pis_4(S^3)\simeq\pis_4(S^2)\simeq\mZ/2$ \cite[Theorems~15.1,\,15.2]{hu}.
 Therefore, the exact sequence for $\pis_4$ arising from \eqref{e31} shows that $\pis_4(M^3(q))=0$
 for $q$ odd and, for $q=2^s$, there is an exact sequence
 \[
  0\to \mZ/2\to \pis_4(M^3_s) \to \mZ/2\to 0.
 \]
 Note that $\pis_4(M^3_1)\simeq\pi_6(M^5_1)$, so \cite[Lemma~10.2]{hu} implies that it embeds
 into $\pi_6(S^3)\simeq\mZ/12$ \cite[Theorem~16.1]{hu}. Hence, $\pis_4(M^3_1)\simeq\mZ/4$.
 For $r>1$ consider the commutative diagram of triangles
 \begin{equation}\label{e33}
  \begin{CD}
   S^2 @>2>> S^2 @>>> M^3_1 @>>> S^3 \\
  @V1VV	@V2^{r-1}VV	@VVV	@V1VV	\\	
   S^2 @>2^r>> S^2 @>>> M^3_s@>>> S^3 , 
 \end{CD} 
 \end{equation}
 It induces the commutative diagram with exact rows
 \[
 \begin{CD}
  0 @>>> \mZ/2 @>>> \pis_4(M^3_1) @>>> \mZ/2 @>>> 0 \\
  && @V0VV	@VVV 	@VV1V \\
  0 @>>> \mZ/2 @>>> \pis_4(M^3_s) @>>> \mZ/2 @>>> 0 ,
 \end{CD} 
 \]
 which shows that the second row is the pushdown of the first one along zero map,
 hence, it splits.
 \end{proof}

  \begin{prop}\label{33}
 \[
  \hos(M^d_s,M^d_r) \simeq  \begin{cases}
  \mZ/4 &\text{ \em if } r=s=1,\\
  \mZ/2\+\mZ/2^m &\text{\em otherwise, where } m=\min(r,s).
 \end{cases}
 \]
 \end{prop}
 \begin{proof}
 Obviously, we may suppose that $m=3$. Since $\pis_4(M^3_1)\simeq\mZ/4$ is
 a module over the ring $\hos(M^3_1,M^3_1)$,  $2\hos(M^3_1,M^3_1)\ne0$,
 hence, $\hos(M^3_1,M^3_1)\simeq\mZ/4$. On the other hand, applying the functor
 $\hos(\,\_\,,M^3_1)$ to the diagram \eqref{e33} with $s>1$, we get a commutative
 diagram with exact rows 
 \[
   \begin{CD}
  0 @>>> \mZ/2 @>>> \hos(M^3_1,M^3_1) @>>> \mZ/2 @>>> 0 \\
    && @A1AA	@AAA 	@AA0A \\
  0 @>>> \mZ/2 @>>> \hos(M^3_s,M^3_1) @>>> \mZ/2 @>>> 0.
 \end{CD} 
 \]
 Thus its second row is the pull-back of the first one along the zero map, hence, it splits.
 The dual consideration shows that the sequence \eqref{e32} for $r>1$
 can be obtained as a pushdown of the sequence for $r=1$, hence, it splits too.
 \end{proof}
  
 Note that the latter decomposition in this statement is that of \emph{groups}. Taking into account
 the multiplication, it is convenient to present morphisms $M^d_s\to M^d_r$ as triangular matrices
 $\mtr{a&b\\0&c}$, with $a\in\mZ/2^r,\,b\in\mZ/2,\,c\in\mZ/2^s$,
 $2^{s-m}a\equiv 2^{r-m}c \mod 2^\mu$, where $m=\min(s,r),\,\mu=\max(s,r)$.
 The product of morphisms correspond then to the usual product of matrices, while the sum of
 morphisms correspond to the usual sum of matrices, with the only exception, when $s=r=1$:
 then we must add matrices as follows: 
 \[
  \mtr{a&b\\0&a}+\mtr{a'&b'\\0&c'}=\mtr{a+a'&b+b'+aa'\\0&a+a'}. 
 \]

\medskip
 Let now $n=3,m=1$, then $\cS_3=\cA_{3,1}\dg\cB_{3,1}$, where $\cA_{3,1}=S^3\cS_1$
 and $\cB_{3,1}=S\cS_2$. Hence, polyhedra from $\cA_{3,1}$ are just bouquets of spheres $S^4$,
 while those from $\cB$ are bouquets of spheres $S^4,S^3$ and Moore spaces $M^4(q)$. 
 For convenience, we set $M^4_0=S^4$ and $M^4_\8=S^3$ and order the set of indices
 by the rule $1<2<\dots<\8<0$.
 As we have seen, $\hos(S^4,M^4(q))=0$ for $q$ odd, $\hos(S^4,M^4_r)=H_r\simeq\mZ/2$ for
 $r\ne 0$ and $\hos(S^4,S^4)=H_0\simeq\mZ$. Therefore, a map $a:A\to B$, where $A\in\cA$,
 $B\in\cB$ can be presented as a block matrix
 \begin{equation}\label{e34}
    a=\mtr{a_0\\ a_\8\\ \vdots\\ a_2\\ a_1} ,
 \end{equation}
  where $a_s$ are matrices over $H_s$. One easily sees that if $\eta_s$ is a generator of
 $H_s$ and $\be_{rs}:M^4_s\to M^4_r$, then $\be_{rs}\eta_s=0$ if $r>s$, while for $r\le s$
 the map $\be_{rs}$ can be so chosen that $\be_{rs}\eta_s=\eta_r$. Set 
 \[
  H_{rs}= \begin{cases}
  0 &\text{ if } r>s,\\
  \mZ/2 &\text{ if } 0\ne r\le s,\\
  \mZ &\text{ if } s=r=0.
 \end{cases} 
 \]
 Therefor two matrices $a,a'$ of the form \eqref{e34} define isomorphic objects from $\bim(\cS_{3,1})$
 \iff there is an invertible integral matrix $\al$ and an invertible block matrix $\be=\be_{rs}$, where
 $\be_{rs}$ is a matrix over $H_{rs}$, such that $a'=\be a\al^{-1}$. Then simple considerations
 show that every object from $\bim(\cS_{3,1})$ decomposes into a direct sum of objects given by the
 $1\xx1$ matrices $q\in H_0$, $\eta_r\in H_r,\ r\ne0$ and $\mtr{2^s\\ \eta_r}\in H_0\+H_r,\
 r\ne0, s>0$. The first case correspond to the Moore space $M^5(q)$, while the second and the third cases
 define new polyhedra, respectively, $C^5(\eta)$, $C^5(2^r\eta)$, $C^5(\eta2^s)$ and $C^5(2^r\eta2^s)$,
 given by the gluing diagrams
 
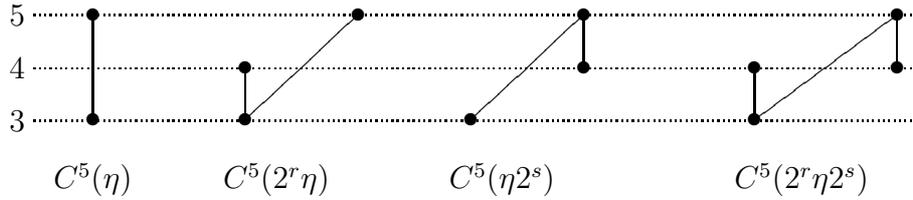
\[
  \xymatrix@R=1ex@C=.7ex{
	 {5} \ar@{.} [rrrrrrrrr] & \cel \ar@{-} [dd] && \cel \ar@{-}[ddl] & & 
	\cel \ar@{-}[ddl]\ar@{-}[d]  &&& \cel \ar@{-}[ddl] \ar@{-}[d] & {} \\
        {4} \ar@{.} [rrrrrrrrr] && \cel \ar@{-}[d] &&& \cel 
	 && \cel\ar@{-}[d] &\cel&{} \\
        {3} \ar@{.} [rrrrrrrrr] & \cel & \cel && \cel &&& \cel&&  {} 	\\
	& {C^5(\eta)}  & {\qquad C^5(2^r\eta)} && {\qquad C^5(\eta2^s)} &&& {\qquad\quad C^5(2^r\eta 2^s)} 	}   
 \]
 (The words in brackets show the corresponding gluings.) 
 
 To find endomorphisms of these atoms, note that there are triangles
 \begin{align}
  & S^3\vee S^4 \xarr{(2^r\ \eta)} S^3 \to C^5(2^r\eta) \to S^4\vee S^5, \label{e35}\\
  & S^4\xarr{\mtr{\eta\\ 2^s}} S^3\vee S^4\to C^5(\eta 2^s) \to S^5, \label{e36} \\
  & S^3\vee S^4 \xarr{\mtr{2^r&\eta\\0&2^s}} \to C^5(2^r\eta2^s)\to S^3\vee S^4. 
 \label{e37}
 \end{align} 
 By Theorem~\ref{11}, $\es(C^5(2^r\eta))$, up to an ideal $I$ such that $I^2=0$, is isomorphic to
 the endomorphism ring of the map $f=(2^r\ \eta)$ in the category $\bim(\cS)/\cJ$. An endomorphism
 of $f$ in $\bim(\cS)$ is a pair $(\al,\be)$, where $\al=\mtr{a&b\eta\\0&c}$ $(a,c\in\mZ,\,b\in\mZ/2)$,
 $\be\in\mZ$, such that $\be f=f\al$, i.e. $\be=a\equiv c\mod2$. Moreover, one easily sees that
 $\cJ$ consists of the pairs with the first component $\mtr{2^rx & x\eta\\ 0&0}$, whence
 $\es(C^5(2^r\eta))/I^2$ is isomorphic to the subring of $\mZ/2^{r+1}\+\mZ$ consisting of all pairs
 $(a,c)$ with $a\equiv c\mod2$. This ring has no nontrivial idempotent, hence, $C^5(2^r\eta)$
 is indeed indecomposable, hence, an atom. Moreover, using the triangle \eqref{e35}, one can see
 that $I\simeq\mZ/2$ and $\es(C^5(2^r\eta))$ is isomorphic to the ring of triangular matrices
 $\mtr{a&b\\0&c}$, where $a\in\mZ/2^{r+1},\,b\in\mZ/2,\,c\in\mZ,$ $a\equiv c\mod2$.
 The same result for $\es(C^5(\eta2^r))$ follows from the triangle \eqref{e36}. Finally, one gets
 from the triangle \eqref{e37} that $\es(C^5(2^r\eta2^s))$ is isomorphic to the ring of triangular
 matrices $\mtr{a&b\\0&c}$, where $a\in\mZ/2^r,\,b\in\mZ/2,\,c\in\mZ/2^s$,
 $a\equiv c\mod2$. Therefore these polyhedra are also atoms. They are called \emph{Chang atoms}
 Moreover, the last ring is local, thus the multiplicity of $C^5(2^r\eta2^s)$ (as well as of any its shift)
 in a decomposition of a polyhedron into a bouquet of indecomposables is the same for all such
 decompositions. Note that the same is true for suspended atoms $M^d(q)$. On the other hand,
 the triangles \eqref{e35} and \eqref{e36} imply that $\rH_3(C^5(2^r\eta)\simeq\rH_4(C^5(\eta2^r)
 \simeq\mZ/2^r$, while other homologies of these spaces are zero. Altogether, it gives
 the following description of the category $\cS_3$.

 \begin{theorem}[Whitehead--Chang, \cite{wh,ch}]\label{34}
  Any polyhedron from $\cS_3$ uniquely (up to permutation of summands) decomposes into a bouquet
 of spheres $S^3,S^4,S^5$, suspended Moore atoms $M^4(q),M^5(q)$ and Chang atoms
 $C^5(\eta)$, $C^5(2^r\eta)$, $C^5(\eta2^s)$ and $C^5(2^r\eta2^s)$.
 \end{theorem}

 Using terms from the representation theory, one can say that the categories $\cS_n,\ n\le3$, are
 \emph{discrete} (or \emph{essentially finite}). In this context it means that there are only finitely
 many isomorphism classes of polyhedra in $\cS_n$ with a prescribed exponent of the torsion part of
 homologies. (So it looks similar to the description of finitely generated abelian groups.)

\section{Tame case: Baues--Hennes Theorem}
 \label{s4}

 We study now the category $\cS_4$. By Theorem~\ref{22},
 $\cS_4=\cA\dg\cB$, where $\cA=S^3\cS_2,\ \cB=S^2\cS_2$.  By
 Theorem~\ref{31}, every polyhedron from $\cA$ (from $\cB$) is a
 bouquet of spheres $S^5,S^6$ and Moore atoms $M^6(q)$ (respectively,
 $S^4,S^5$ and $M^5(q)$).  We have already calculated morphisms
 between indecomposables in $\cS_2$. Just in the same way one
 calculates morphisms from the objects of $S^3\cS_2$ and those of
 $S^2\cS_2$.  We omit the details, which are standard; the result is
 presented in Table~1.
 \begin{table}[ht]
 \caption{}\vspace*{-1.5em}
 \[
   \begin{array}{|c|cccc|}
  \hline 
 \phan & S^5 & S^6 & M^6_1 & M^6_s\ (s>1) \\
 \hline 
   S^4 &\phan  \mZ/2 & \mZ/2 & \mZ/4 & \mZ/2\+\mZ/2 \\
  S^5 &\phan  \mZ & \mZ/2 & \mZ/2 & \mZ/2 \\
   M^5_1 &\phan  \mZ/2 & \mZ/4 & \mZ/2\+\mZ/2 & \mZ/4\+\mZ/2 \\
   M^5_r\ (r>1) &\phan  \mZ/2 &\ \mZ/2\+\mZ/2\ &\ \mZ/2\+\mZ/4\ &\
 \mZ/2\+\mZ/2\+\mZ/2  \\ 
 \hline
 \end{array}
 \]
 \end{table}
 Actually, the groups $\hos(M^6_s,M^5_r)$ can be naturally considered
 as the groups of upper triangular matrices $\mtr{a&b\\0&c}$ over
 $\mZ/2$ with $b=0$ if $s=r=1$.  Again the sum of morphisms correspond
 to the usual sum of matrices, with the exceptions for $s>1,r=1$ and
 $s=1,r>1$, when the sum of matrices must be twisted as follows:
 \begin{align*}
  \mtr{a&b\\0&c}+\mtr{a'&b'\\0&c'}&=\mtr{a+a'&b+b'+aa'\\0&c+c'} 
   \text{ if } s>1,r=1,\\
 \mtr{a&b\\0&c}+\mtr{a'&b'\\0&c'}&=\mtr{a+a'&b+b'+cc''\\0&c+c'} 
   \text{ if } s=1,r>1.
 \end{align*}
 The multiplication of elements from $\hos(M^6_s,M^5_r)$ by morphisms
 between objects from $\cA$ and $\cB$ (also presented by triangular
 matrices as in Section~\ref{s3}) correspond to the usual product of
 matrices. Therefore, a morphism $A\to B$ can be naturally considered
 as a block matrix presented in Table~2.
 \begin{table}[htp]
 \caption{}\vspace*{-1.5em}
  \[
  \left\uparrow \underrightarrow{\hspace*{.5ex}\left[
          \begin{array}{l|ccccc|ccccccc}
	    \phantom{\Big|}
 	 & \nc1 & \nc2 & \nc3 &\dots &  &&  &\dots & \nc3 &\nc2 &\nc1
 	 \\
	 \hline
	\hskip-4pt \phantom{\Big|}
     	\nn1&2 & 2 & 2 & \dots & 2 && 2 &\dots & 2 & 2 & 0 &\\
	\nn2 &2& 2 & 2 & \dots & 2 && 2 &\dots & 2 & 2 &2 &\\
	 \nn3&2& 2 & 2 & \dots & 2 && 2 &\dots & 2 & 2 & 2\\
	&\vdots &\vdots &\vdots &\ddots &\vdots & &\vdots &\ddots &
	\vdots &\vdots &\vdots &\\
	&2& 2 & 2 & \dots & 2 && 2 &\dots & 2 & 2 & 2&\\
 \hline 
	\phantom{\Big|} &0 & 0 & 0 & \dots & \8 && 2 &\dots & 2 & 2 & 2 &\\
	&\vdots &\vdots &\vdots &\ddots &\vdots & &\vdots &\ddots &
	\vdots &\vdots &\vdots &\\
	\nn{3} &0 & 0 & 0 & \dots & 0 && 2 &\dots & 2 & 2 & 2 &\\
	 \nn{2}&0 & 0 & 0 & \dots & 0 && 2 &\dots & 2 & 2 & 2 &\\
	\nn{1} &0 & 0 & 0 & \dots & 0 && 2&\dots & 2& 2& 2 &
 \end{array} \right] \vspace*{.7ex} } \right. 
  \]
 \end{table}
  In this table a symbol $2$ ($\8$) shows that the corresponding block
 has values from $\mZ/2$ (respectively, from $\mZ$). Zeros show that
 the corresponding block is always zero.  Arrows on the left and
 below symbolize the action of morphisms between the objects from
 $\cA$ and $\cB$ respectively. The labels $\nn1,\nn2,\ldots$
 (or $\nc1,\nc2,\ldots$) show that the corresponding horizontal
 (respectively, vertical) stripes are of the same size and we must use
 the same elementary transformations in both of them. These stripes
 correspond to $M^d_r$ with the same $d$ and $r$. Note that there
 are 2 horizontal and 2 vertical stripes without such labels. They
 correspond to spheres $S^d$.

 This matrix problem is a slight variation of a well-known one,
 namely, representations of \emph{bunches of chains} (see \cite{bon}
 or \cite[Appendix B]{bd}). It implies a description of indecomposable
 objects in the category $\bim(\cS_{4,2})$, hence, in $\cS_4$. We call
 them \emph{strings} and \emph{bands}, as it is usual in the
 representation theory of algebras. Not providing details (see
 \cite{d1}), we just present the corresponding attachment diagrams
 (Table~3). It is convenient to distinguish two types of strings:
 \emph{usual} and \emph{decorated}; I hope that the pictures show the
 difference. ``Decorations'' (one for each string) are shown with
 double lines. We omit integers precising the degrees of ``vertical''
 attachments, as well as one precising the ``long'' attachment in a
 decorated strings of the first kind; they can be arbitrary and differ
 for different attachments. Certainly, each diagram is actually
 finite: it starts at any place on the left and stops at any place on
 the right.
 \begin{table}[htp]
   \caption{}\vspace*{-1.3em}
  \medskip
   \begin{center}
       \sf usual strings
   \end{center}
\medskip
 \[
 \xymatrix@R=.1ex@C=3ex@!R{
	 {7} \ar@{.} [rrrrrrrrrrrrrr] &&&& \cel \ar@{-}[dd] &&&& \cel \ar@{-}[dd]
 &&&& \cel \ar@{-}[dd] && {} \\
  {}\\
	{6} \ar@{.} [rrrrrrrrrrrrrr] &&&& \cel \ar@{-}[ddddrr] &&&& \cel \ar@{-}[ddddrr] 
 &&&& \cel \ar@{-}[ddr] && {}  \\
	& {\cdots\quad } \\
	{5} \ar@{.} [rrrrrrrrrrrrrr] & *=0{}& \cel \ar@{-}[uuuurr] \ar@{-}[dd] &&&&
 \cel \ar@{-}[uuuurr] \ar@{-}[dd] &&&&  \cel \ar@{-}[uuuurr] \ar@{-}[dd] &&& *=0{} & {} \\
	&&&&&&&&&&&&& {\quad \cdots}  \\
	{4}  \ar@{.} [rrrrrrrrrrrrrr] && \cel@\ar@{-}[uul] &&&& \cel &&&& \cel &&&& {}
 } 
 \]
 \medskip
 \begin{center}
     \sf decorated strings
 \end{center}
\medskip
 \[
 \hspace*{-1cm}
    \xymatrix@R=.1ex@C=3ex@!R{
	 {7} \ar@{.} [rrrrrrrrrrrrrrrr] && \cel \ar@{-}[dd] \ar@{-}[ddddrr] &&&&
 \cel \ar@{-}[dd] \ar@{-}[ddddrr] &&&& \cel \ar@{-}[dd] &&&& \cel \ar@{-}[dd] && {} \\
 {} \\
	{6} \ar@{.} [rrrrrrrrrrrrrrrr] && \cel \ar@{-}[ddl] &&&& \cel \ar@{-}[ddddll] 
 &&&& \cel \ar@{-}[ddddrr] &&&& \cel \ar@{-}[ddr] && {}  \\
	{} \\
	{5} \ar@{.} [rrrrrrrrrrrrrrrr] & *=0{}&&& \cel \ar@{-}[dd] &&&&
 \cel \ar@{-}[dd] &&&&  \cel \ar@{-}[uuuurr] \ar@{-}[dd] &&& *=0{} & {} \\
	& {\cdots\quad} &&&&&&&&&&&&&& {\quad \cdots}  \\
	{4}  \ar@{.} [rrrrrrrrrrrrrrrr] &&&& \cel &&&& \cel
    \ar@{=}[uuuuuurr] &&&& \cel &&&& {} 
 }  
 \] 

\smallskip
\begin{flushleft}
  and
\end{flushleft}
\smallskip
 \[ \hspace*{-1cm}
      \xymatrix@R=.1ex@C=3ex@!R{
	 {7} \ar@{.} [rrrrrrrrrrrrrrrr] &&&& \cel \ar@{-}[dd] \ar@{-}[ddddrr] &&&&&&
  \cel \ar@{-}[dd] \ar@{=}[ddddll] &&&& \cel \ar@{-}[dd] \ar@{-}[ddddll]  & & {} \\
  &{\cdots\quad} \\
	{6} \ar@{.} [rrrrrrrrrrrrrrrr] &&&& \cel \ar@{-}[ddddll] &&&&
 \cel \ar@{=}[ddddll] \ar@{=}[dd] &&  \cel \ar@{-}[ddddrr]  &&&& \cel \ar@{-}[ddr] && {} \\
  {} \\
	{5} \ar@{.} [rrrrrrrrrrrrrrrr] && \cel \ar@{-}[uul]
	 \ar@{-}[dd] &&&& \cel \ar@{-}[dd]  
 && \cel  &&&& \cel \ar@{-}[dd] &&& *=0{}& {}  \\
	& &&&&&&&&&&&&&& {\quad \cdots}  \\
	{4}  \ar@{.} [rrrrrrrrrrrrrrrr] && \cel &&&& \cel  &&&&&& \cel &&&& {}
 }   
 \]
 \medskip
 \begin{center}
   \sf bands
 \end{center}
 \medskip
 \[
 \xymatrix@R=.1ex@C=3ex@!R{
	 {7} \ar@{.} [rrrrrrrrrrrrrr] &&& \cels \ar@{-}[dd] &&&& \cels \ar@{-}[dd]
 &&&&&& \cels \ar@{-}[dd] && {} \\
  {}\\
	{6} \ar@{.} [rrrrrrrrrrrrrr] &&& \cels \ar@{-}[ddddrr] &&&& \cels \ar@{-}[dr] 
 &&&&&& \cels  && {}  \\
	&&&&&&&& *=0{}& {\hspace*{-1em}\cdots}&&*=0{}  \\
	{5} \ar@{.} [rrrrrrrrrrrrrr] & \cels \ar@{-}[uuuurr] \ar@{-}[dd] &&&&
 \cels \ar@{-}[uuuurr] \ar@{-}[dd] &&&&&&  \cels \ar@{-}[uuuurr]
 \ar@{-}[dd] &&& *=0{} & {} \\ 
	&&&&&&&&&&&&&   \\ 
	{4}  \ar@{.} [rrrrrrrrrrrrrr] &
 \cels\ar@{~}[uuuurrrrrrrrrrrr]|(.8){\mbox{\boldmath$\Phi$}} &&&& \cels &&&&&& \cels
 \ar@{-}[uull] &&&& {} 
 }
 \]
 \end{table}

 Multiple bullets in the case of bands symbolize not a unique
 cell but several (say $m$) copies of it (the same for each
 ball). All attachments except the one marked by the wavy line are
 ``natural'': the first copy of an upper cell is attached to the first
 copy of a lower one, the second to the second, etc. The attachment
 marked by the wavy line is ``twisted'' by an invertible Frobenius
 matrix $\Phi$ of size $m\xx m$ over the field $\mZ/2$ with the
 characteristic polynomial $f(x)$, which must be a power of an irreducible
 one and such that $f(0)\ne0$. For instance, if $f(x)=x^3+x+1$,
 i.e. $m=3$ and $\Phi=${\normalsize$\mtr{0&0&1\\1&0&0\\0&1&1}$}, this
 attachment is:

 
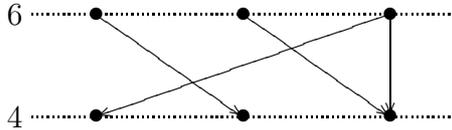
\begin{center}
\[
 \xymatrix{
 {6} \ar@{.} [rrrrrr]    & \cel \ar[drr] && \cel \ar[drr] && \cel \ar[dllll] \ar[d] & \\
 {4} \ar@{.} [rrrrrr]    & \cel && \cel && \cel &   
 }
\]
 \end{center}

One can check that all strings and bands are indecomposable and
pairwise non-isomorphic.  Note also that all atoms from $\cS_4$ are
$p$-primary ($2$-primary, except Moore atoms $M^d(p^r)$ with odd $p$,
which are $p$-primary). Therefore, we have the uniqueness of
decomposition of spaces from $\cS_4$ into bouquets of suspended
atoms. So we get the following result. We call strings and bands
\emph{Baues atoms}.

\begin{theorem}[Baues--Hennes \cite{bh}]
 Any polyhedron from $\cS_4$ decomposes uniquely into a bouquet of
 spheres, suspended Moore atoms, suspended Chang atoms and Baues atoms.
\end{theorem}

In Section~\ref{s7} we shall see that actually $\cS_4$ is the last case
where a ``good'' description of polyhedra is possible. Starting from
$\cS_5$ this problem becomes \emph{wild}.

 \section{Torsion free polyhedra. Finite case}
 \label{s5}

 Consider now torsion free case. Note that if all $\rH_k(X)$ are
 torsion free, the attachment diagram cannot contain ``Moore
 fragments'' 
 \[
  \xymatrix@=1ex{ \ar@{.}[rr] & \cel \ar@{-}[d] & \\
    \ar@{.}[rr] & \cel  & }
 \]
 In particular, among the atoms from Sections~\ref{s3} and \ref{s4}
 only Chang atom $C^5=C^5(\eta)$ and the \emph{double Chang atom}
 $C^7_2=C^7(\eta^2)$ with the attachment diagram
 
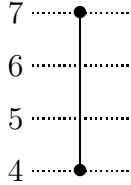
\[
  \xymatrix@R=1ex@C=.7ex{7 \ar@{.}[rrrr] && \cel \ar@{-}[ddd] && \\
    6 \ar@{.}[rrrr] &&&& \\
    5 \ar@{.}[rrrr] &&&& \\
    4 \ar@{.}[rrrr] && \cel  && }
 \]
 are torsion free. Therefore, if we set in Theorem~\ref{24} $n=5,m=4$,
 the category $\cA^0$ consists of bouquets of spheres $S^8$ and the
 category $\cB^0$ consists of bouquets of spheres $S^d\ (5\le d\le
 8)$, suspended Chang atoms $C^7,C^8$ and suspended double Chang atoms
 $C^8_2$. Obviously, $\cS^0(S^8,S^8)=0$. Easy calculation give the
 following values of the groups $\Ga(B)=\cS^0(S^8,B)$ for atoms $B$ from
 $\cB^0$: 

 \[
 \begin{array}{|c|cccccc|}
 \hline
 \phan B\, &\ S^5\ &\ S^6\ &\ S^7\ &\ C^7\ &\ C^8\ &\ C_2^8 \\
 \hline 
 \phan\Ga & \mZ/24 & \mZ/2 & \mZ/2 & \mZ/12 & 0 &\mZ/12 \\
 \hline 
 \end{array}
 \]

\smallskip
 \noindent
 Morphisms of these spaces induce monomorphisms
 $\Ga(S^7)\to\Ga(C^7)\to\Ga(S^5)$ and $\Ga(S^6)\to\Ga(S^5)$,
 epimorphism  $\Ga(S^5)\to\Ga(C_2^8)$, and isomorphisms
 $\Ga(S^7)\to\Ga(S^6)$ and $\Ga(S^7)\to\Ga(S^6)$. Thus, an object from
 $\bim(\cS^0_{5,4})$ can be presented by a block matrix as in
 Table~4.
\begin{table}[htp]
 \caption{}\vspace*{-1.5em}
   \[
  \xymatrix@R=0em@C=1.5ex {*=0{} \ar[d] & *+<2.5em,2em>[F]{\mZ/24}& *=0{} \\
              *=0{}\ar[d]&  *+<2.5em,2em>[F]{\mZ/12} & *=0{} \\
                *=0{}       &  *+<2.5em,2em>[F]{\mZ/12} & *=0{}\ar[uu]|2 \\
                *=0{}\ar[d]& *+<2.5em,2em>[F]{\,\mZ/2\,} & *=0{}\ar[u]|6  \\
               *=0{}       &   *+<2.5em,2em>[F]{\,\mZ/2\,} & *=0{} \ar@/_/@<-1ex>[uuuu]|{12}
	     }
 \]
\end{table}
 Here inside each blocks we have written the groups,
 wherefrom the coefficients of this block are. The arrows show the
 allowed transformation between blocks. An integer $k$ in the arrows
 point out that, when we perform this transformation, the row must be
 multiplied by $k$. (No integer means that $k=1$.) For instance, we
 can add the rows of the third matrix multiplied by $2$ to the rows of
 the first one. Certainly, compositions of these transformations are
 also allowed. Thus, for instance,  we can add the rows of the third
 matrix multiplied by $2$ to the rows of the second one too. The
 arising matrix problem is rather simple. It is of finite type, and
 Table~6 shows the attachment diagrams of the corresponding atoms
 from $\cT_5$. We call them \emph{A-atoms of the 1st kind}.
\begin{table}[!htp]
  \caption{}\vspace*{-1.5em}
 \begin{align*}
   &\xymatrix@R=.4em@C=1em{
     9 \ar@{.}[rrrrrrrrrr] & \cel \lin[dddd] &&  \cel \lin[ddddl] &&
     \cel \lin[ddddl] && \cel \lin[ddddl] \lin[dd] && \cel
     \lin[ddddl] \lin[ddd] & \\
     8 \ar@{.}[rrrrrrrrrr] &&&& \cel \lin[ddd] &&&&&& \\
     7 \ar@{.}[rrrrrrrrrr] && \cel \lin[dd] &&&&& \cel &&& \\
     6 \ar@{.}[rrrrrrrrrr] &&&&&&&&& \cel & \\
     5 \ar@{.}[rrrrrrrrrr] & \cel & \cel && \cel && \cel && \cel && \\
     & A(v) & A(\eta v)&& A(\eta^2v) && A(v\eta) && A(v\eta^2)
     } \\
   \\
    & \xymatrix@R=.4em@C=1em{
     9 \ar@{.}[rrrrrrrrr] && \cel \lin[dd] \lin[ddddl] &&  \cel
     \lin[ddd] \lin[ddddl] && \cel \lin[ddddl] \lin[dd] && \cel
     \lin[ddddl] \lin[ddd]] & \\
     8 \ar@{.}[rrrrrrrrr] & \cel \lin[ddd] && \cel \lin[ddd]
     &&&&&& \\
     7 \ar@{.}[rrrrrrrrr] && \cel  &&& \cel  \lin[dd] & \cel &\cel
     \lin[dd] && \\ 
     6 \ar@{.}[rrrrrrrrr] &&&& \cel &&&& \cel & \\
     5 \ar@{.}[rrrrrrrrr] & \cel && \cel && \cel && \cel &&  \\
     & A(\eta^2 v\eta) && A(\eta^2 v\eta^2)&& A(\eta v\eta)
     && A(\eta v\eta^2) &&
     } \\
  \end{align*}
  \end{table}
 The integer $v$ show, which multiple of the generator of the group
   $\pis_8(S^5)\simeq\mZ/24$ is used for the ``long'' attachment. Actually,
  $ 1\le v\le12$ in the case of $A(v)$, $1\le v\le 3$ in the case
 of $A(\eta v\eta)$, $1\le v\le 6$ in all other cases.   
 
   So we have got a description of polyhedra from $\cT_5$.

 \begin{theorem}[Baues--Drozd \cite{bd1}]
   Every polyhedron from $\cT_5$ is a bouquet of spheres, suspended
   Chang and double Chang atoms, and the A-atoms of the first kind.
 \end{theorem}

 Note that this time the decomposition is not unique; even the
 cancellation law does not hold. For instance,
 $A(3)\+S^5\simeq A(9)\+S^5$ \cite{bd1,d1}; see ibidem more on
 decomposition laws.

 Analogous is the case of $\cT_6$, when we take $m=4$. We omit
 details, just schematically presenting the arising matrix problem in Table~6.
 The dashed line from the 4th to the 6th level show the transformation that
 only acts on the left-hand column (on $\mZ/2$ components).
 \begin{table}[!htp]
   \caption{}\vspace*{-1.5em}
   \[
   \xymatrix@R=0em@C=0em {*=0{} \ar@<-1em>[d] &
   *+<2.5em,2em>[F]{\mZ/24}&*+<3.4em,2.25em>[F]{\ 0\ } & *=0{} \\ 
              *=0{}\ar@<-1em>[d]&
   *+<2.5em,2em>[F]{\mZ/12}&*+<3.35em,2.45em>[F]{\ 0\ } & *=0{} \\ 
                *=0{}       &
   *+<2.5em,2em>[F]{\mZ/12}&*+<3.35em,2.45em>[F]{\ 0\ } &
   *=0{}\ar@<-1em>[uu]|2 \\ 
                *=0{}\ar@<-1em>[d] \ar@{-->}@/_/@<-1.5em>[dd]&
   *+<2.6em,2em>[F]{\,\mZ/2\,}&*+<2.6em,2em>[F]{\,\mZ/2\,} &
   *=0{}\ar@<-1em>[u]|6  \ar@<1em>@{|->}[ddd]|6\\ 
               *=0{}\ar@<-1em>[d]|{12}       &
   *+<3.4em,2.2em>[F]{\ 0\ }&*+<2.6em,2em>[F]{\,\mZ/2\,} & *=0{}
         \\
                  *=0{}  \ar@<-1em>[d]      &
   *+<2.6em,2em>[F]{\,\mZ/2\,}&*+<2.5em,2em>[F]{\mZ/24} & *=0{}
   \ar@<-2.5em>[uuuuu]|{12} \\
                *=0{}       &
   *+<3.4em,2.2em>[F]{\ 0\ }&*+<2.5em,2em>[F]{\mZ/12} & *=0{}
   	\ar@<-2.5em>[u]|2     }
  \]
 \end{table}
  The resulting list of atoms (their attachment diagrams) see in the Table~7.
 \def\cell{\cel \lin[ddddl]}
 \begin{table}[!htp]
  \caption{}\vspace*{-1.4em}
 \begin{align*} 
 & \hspace*{-2em}
  \xymatrix@R=.4em@C=1.5em{
  11 \ar@{.}[rrrrrrrrrrrrrrrr] &&& \cell &&& \cell \lin[dd] &&&
   \cell \lin[ddd] &&& \cell \lin[dd] &&& \cell & \\
  10 \ar@{.}[rrrrrrrrrrrrrrrr] && \cell \lin[ddd] &&& \cell \lin[ddd]
 &&& \cell \lin[ddd] &&& \cell \lin[ddd] &&& \cell \lin[ddd] && \\
   9 \ar@{.}[rrrrrrrrrrrrrrrr] &&&&&& \cel & \cel \lin[ddd] &&&&&
	\cel &&&& \\
   8 \ar@{.}[rrrrrrrrrrrrrrrr] &&&& \cel \lin[dd] &&&&& \cel &&&&
	\cel \lin[dd] &&&& \\
   7 \ar@{.}[rrrrrrrrrrrrrrrr] && \cel &&& \cel &&& \cel &&& \cel &&&
	\cel && \\
   6 \ar@{.}[rrrrrrrrrrrrrrrr] & \cel &&& \cel &&& \cel &&& \cel &&&
	\cel &&&   
	  \\  \scriptstyle
      &&  *=0{A(v\eta^2w)} &&& *=0{A(\eta v\eta^2w\eta)} &&&
   *=0{A(\eta^2v\eta^2w\eta^2)} &&& *=0{A(v\eta^2w\eta)} &&&
   *=0{A(\eta v\eta^2 w)} &&
	}\\ \\
 & \xymatrix@R=.4em@C=1.5em{
 11 \ar@{.}[rrrrrrrrrrrrr] &&& \cell \lin[ddd] &&& \cell &&& \cell
  \lin[dd] &&& \cell \lin[ddd] & \\
 10 \ar@{.}[rrrrrrrrrrrrr] && \cell \lin[ddd] &&& \cell \lin[ddd] &&&
  \cell \lin[ddd] &&& \cell \lin[ddd] && \\
  9 \ar@{.}[rrrrrrrrrrrrr] &&&& \cel \lin[ddd]  &&& \cel \lin[ddd] &&
  \cel &&&& \\
  8 \ar@{.}[rrrrrrrrrrrrr] &&& \cel &&&&&&& \cel \lin[dd] && \cel & \\
  7 \ar@{.}[rrrrrrrrrrrrr] && \cel &&& \cel &&& \cel &&& \cel && \\
  6 \ar@{.}[rrrrrrrrrrrrr] & \cel &&& \cel &&& \cel &&& \cel &&& \\
  \scriptstyle
  && *=0{A(v\eta^2w\eta^2)} &&& *=0{A(\eta^2v\eta^2w)} &&&
  *=0{A(\eta^2v\eta^2w\eta)} &&& *=0{A(\eta v\eta^2w\eta^2)} && 
   }
 \end{align*}
 \end{table}
  We call them \emph{A-atoms of the second kind}. The integers $v$ and
 $w$ show, as above, the multiple of generator, respectively, of
 $\pis_9(S^6)$ and $\pis_{10}(S^7)$ used for the corresponding
 attachments. In all cases $v,w\in\set{1,2,3,4,5,6}$.

    So we have got a description of polyhedra from $\cT_6$.

  \begin{theorem}[Baues--Drozd \cite{bd4}]
   Every polyhedron from $\cT_6$ is a bouquet of spheres, suspended
   Chang and double Chang atoms, suspended A-atoms of the first kind
   and A-atoms of the second kind.
 \end{theorem}

 In the next section we shall use the values of $\hos$-groups between
 Chang atoms and spheres. To deal with the Chang atom $C^5$ we apply the
 bifunctor $\hos$ to the cofibration sequence
 	\begin{equation}\label{fib}
 	S^4 \xarr\eta S^3 \to C^5 \to S^5 \xarr\eta S^4.
 \end{equation}
 It gives the commutative diagram with exact rows and columns (we write
  here $(X,Y)$ instead of $\hos(X,Y)$\,)
	\[
	\begin{CD}
	\mZ @>>> \mZ/2 @>>> (C^5,S^4) @>>> 0 @>>> \mZ \\
	@VVV	@VVV 	@VVV	@VVV	@VVV \\
	\mZ/2 @>>> \mZ/2 @>>> (C^5,S^4) @>>> \mZ @>>> \mZ/2 \\
	@VVV	@VVV 	@VVV	@VVV	@VVV \\
	(S^4,C^5) @>>> (S^5,C^5) @>>> (C^5,C^5) @>>>
	 (S^3,C^5) @>>> (S^4,C^5) \\
	@VVV	@VVV 	@VVV	@VVV	@VVV \\
	0 @>>> \mZ @>>> (C^5,S^5) @>>> 0 @>>> 0 \\
	@VVV	@VVV 	@VVV	@VVV	@VVV \\
	\mZ @>>> \mZ/2 @>>> (C^5,S^4) @>>> 0 @>>> \mZ ,
	\end{CD}
\]
 where all maps $\mZ\to\mZ/2$ are surjective and all maps
 $\mZ/2\to\mZ/2$ are bijective. It gives the following values of
 $\hos$-groups:
	\begin{align*}
	\hos(C^5,S^4)&=\hos(S^4,C^5)=0\\
	\hos(S^3,C^5)&=\hos(C^5,S^5)=\mZ,\\
	\hos(C^5,S^3)&=\hos(S^5,C^5)=2\mZ,\\
	\hos(C^5,C^5)&=\mD,
\end{align*}
 where $\mD$ (the ``\emph{dyad}'') is the subrings of $\mZ\xx\mZ$
 consisting of all pairs $(a,b)$ with $a\equiv b\mod2$. 
 
 Similar observations applied to the suspended versions of the sequence
 \eqref{fib} and the cofibration sequence
	\[
	S^6\xarr{\eta^2} S^4 \to C^7_2 \to S^7 \xarr{\eta^2} S^5
\]
 give Table~8 of the values $\hos(X,Y)$ for suspended
 atoms from $\cT_4$.
 \begin{table}[!htp]
 \caption{}\vspace*{-1.4em}
 \[	\hspace*{-1.5ex}
	 \begin{array}{|c|c|cc|cc|c|cc|c|c|}
	 \hline
	 \phan & S^4 & C^7_2: 4 & 7 & C^6: 4& 6 & S^5 & C^7:5 &7 
	 & S^6 & S^7 \\
	 \hline
	 S^4\phan & \mZ & 2\mZ & \mZ/12 & 2\mZ & 0 & \mZ/2 & 0 &
	 \mZ/12 & \mZ/2 & \mZ/24 \\
	 \hline
	 C^7_2: 4\phan & \mZ & \mZ^= & \mZ/12 & 2\mZ & 0 & \mZ/2
	 & 0 & \mZ/12 & 0 & \mZ/12 \\
	 \phh 7\phan & 0 & 0 & \mZ^= & 0 & 0 & 0& 0& \mZ & 0 &2\mZ\\
	 \hline
	 C^6: 4\phan & \mZ & \mZ & \mZ/12 & \mZ^= & 0 & 0 & 0 &
	 \mZ/12 & 0 & 2\mZ \\
	 \phh 6 \phan & 0& 0 & 0 & 0 & \mZ^= & 0 & 0 & 0 & 2\mZ & 0 \\
	 \hline
	 S^5\phan & 0 & 0 & 0 & 0 & 0 & \mZ & 2\mZ & 0 & \mZ/2
	 & \mZ/2 \\
	 \hline
	 C^7: 5\phan & 0 & 0 & 0 & 0 & 0 & \mZ & \mZ^= & 0 &
	 0 & 0 \\
	 \phh 7\phan & 0 &0 & 2\mZ & 0 & 0 & 0 & 0 & \mZ^= &
	 0 & 2\mZ \\  
	 \hline
	 S^6\phan & 0 & 0 & \mZ/2 & 0 & \mZ & 0 & 0 & 0 & 
	 \mZ & \mZ/2 \\
	 \hline
	 S^7\phan & 0& 0 & \mZ & 0 & 0 & 0 & 0 & \mZ & 0 & \mZ
	 \\  \hline
	\end{array}
\]
\end{table}
 In this table the $\hos$-groups for suspended Chang atoms are
 presented in matrix form, emphasizing which components have come
 from the cells of given dimensions. The superscripts $^=$ show
 that the diagonal parts of the corresponding matrices are with 
 entries not from $\mZ\xx\mZ$, but from $\mD$. For instance, 
 $\es(C^7_2)$ is presented as the ring of triangular matrices
 $\mtr{a&b\\0&c}$, where $a,c\in\mZ,\ a\equiv c\mod2,\ b\in\mZ/12$.
 Under such presentation the multiplication of morphisms turns
 into the multiplication of matrices.

\section{Torsion free polyhedra. Tame case}
 \label{s6}

  The category $\cT_7$ is more complicated. To describe it, we use
  Theorem~\ref{24} with $n=7,m=3$. Then $\cA^0$ consists of the
  bouquets of spheres $S^d\ (10\le d\le12)$ and suspended Chang atoms
  $C^{12}$, while $\cB^0$ consists of bouquets of spheres
  $S^d\ (7\le d\le 10)$ suspended Chang atoms $C^9,C^{10}$ and
  suspended double Chang atoms $C^{10}_2$. The calculations similar
  to those of the end of preceding section give Table~9 of the values
  of groups $\cS^0(X,Y)$ for the suspended atoms $X\in\cA^0$ and 
  $Y\in\cB^0$, also presented in matrix form.
  \begin{table}[!htp]
 \caption{}\vspace*{-1.4em}
 \[	\hspace*{-1.5ex}
	 \begin{array}{|c|c|c|c|cc|}
	 \hline
	 \phan & S^{10} & S^{11} & S^{12} & C^{12}:10 &12	 \\
	 \hline
	 S^7\phan & \mZ/24 & 0 & 0 & \mZ/24 & 0  \\
	 \hline
	 C^{10}_2: 7\phan & 0 & 0 & \mZ/2 &  \mZ/24^* & 0 \\
	 \phh 10\phan & 0 & 0 & \mZ/2 & 0 & \mZ/2^* \\
	 \hline
	 C^9: 7\phan & \mZ/12 &  0 & 0 & \mZ/24^* & 0 \\
	 \phh 9 \phan & 0& 0 & \mZ/24 & 0 & \mZ/2^*  \\
	 \hline
	 S^8\phan &  \mZ/2 & \mZ/24 & 0 & 0& 0 \\
	 \hline
	 C^{10}: 8\phan & 0 & \mZ/12 & 0 & 0 & 0 \\
	 \phh 10\phan & 0 & 0 & 0 & 0 & 0  \\  
	 \hline
	 S^9\phan &  \mZ/2 &  \mZ/2 & \mZ/24 & 0 & \mZ/12 \\
	 \hline
	 S^{10}\phan & 0& \mZ/2 & \mZ/2 & 0 & 0
	 \\  \hline
	\end{array}
\]
\end{table}
  The superscripts $^*$ show that in the corresponding groups we
  identify the elements of order $2$. So actually, these values are
  isomorphic to $\mZ/24$, but it is convenient to consider them as
  $(\mZ/24\+\mZ/2)/(12,1)$. Then again the action of morphisms from
  $\cA^0$ and $\cB^0$, as presented in Table~8 (or, rather, its 
  suspended version) turns into the multiplication of matrices.
  Again we obtain a bimodule problem close to that of bunches of
  chains, especially, in its ``\emph{decorated}'' version (see
  \cite{odk}). To present the answer (for details see \cite{d2}), we
  introduce the following notations and definitions.

  \begin{defin}
    \begin{enumerate}
    \item  We consider chains $\dE_k$ and $\dF_k$ ($1\le k\le $):
      \begin{align*}
	\dE_1&=\set{e_1<e_2<e_4}, &\dF_1&=\set{f_4<f_1},\\
	\dE_2&=\set{e_5<e_9}, &\dF_2&=\set{f_3<f_5},\\
	\dE_3&=\set{e_6<e_7}, &\dF_3&=\set{f_2},\\
	\dE_4&=\set{e_3<e_{10}<e'_9<e'_6}, &\dF_4&=\set{f'_1<f'_2<f'_3}.
      \end{align*}
      Actually, the elements $e_i$ ($f_j$) correspond to the rows
  (columns) of Table~9, while the relations $\rr{c}$ correspond to the
  elements of the groups $\cC^0(A,B)$. We need extra elements $e'_i$
  and $f'_j$ since the entries $\mZ/2$ in this table behave in a
  different way than the other ones. 

    We set
      $\dE=\bup_{i=1}^4\dE_i,\,\dF=\bup_{i=1}^4\dF_i,\,\dX=\dE\cup\dF$. 
      $x\approx y$ means that $x$ and $y$ belong to the same set
      $\dE_i$ or $\dF_i$.
    \item  We define symmetric relations $\sim$ and $-$ on $\dX$
    setting $x-y$ if $x\in\dE_i,\,y\in\dF_i$ or vice versa; $e_i\sim
    e_i' (i\in\set{6,9},\ f_j\sim f'_j\ (1\le j\le 3)$. We also define
    the symmetric relations $\rr{c}$, where $c\in\set{1,2,3,4,6}$, setting
    $e_i\rr{c}f_j$ if $e_i-f_j$ and the $(ij)$-th entry in Table~9 is
    $\mZ/m$ with $c\mid m$. We denote by $R$ the set of all relations
    $\{\sim,\rr{c}\}$ and by $v(c)$ the biggest $d$ such that
    $2^d$ divides $c$. 
    \item  We define a \emph{word} as a sequence
	 $ w=   x_1 r_2 x_2 r_3 \dots r_l x_l$
    where $x_i\in\dX,\,r_i\in R$ such that
    \begin{enumerate}
      \item $x_{k-1} r_k x_k$ in $\dX$ for each $1<k\le l$;
      \item if $r_k=\sim$, then $r_{k+1}=\rr{c}$ and vice versa;
      \item if $r_2=\rr{c}$ (respectively, $r_l=\rr{c}$), there is
      no element $y\in\dX$ such that $x_1\sim y$ (respectively,
      $x_l\sim y$);
      \item if $r_k=\rr{c}\,$ with $v(c)=1$, then either $2<k<l$, or
      $k=2,\,x_1=e_1$, or $k=l,\,x_l=x_1$;
      \item if $\fR=\rr{c}$ with $v(c)=2$, then $\fR$ can only occur in
      the following words or their reverse:
      \begin{align*}
	& e_4\sim e_5 \fR f_3 \sim\dots\ \text{ (of any length)},\\
	& e_1\fR f_4\sim f_5,\ e_3\sim e_2\fR f_4\sim f_5,\\
	& \dots\rr{c'}e_4\sim e_5 \fR f_3 \sim\dots\ \text{ (of any
	length)},\\ 
	& e_1\fR f_4\sim f_5\rr{c'}\dots\ \text{ (of any length)},\\
	& e_3\sim e_2\fR f_4\sim f_5\rr{c'}\dots\ \text{ (of any length)},\\
	& e_1\fR f_4\sim f_5\rr{c'}\dots\ \text{ (of any length)},\\
	& e_1\fR f_1\sim f_1',\ e'_6\sim e_6\fR f_2\sim f'_2,\ e'_9\sim
	e_9\fR f_3\sim f'_3, 
      \end{align*}
      where $c'\equiv0\pmod3$;
      \item if $w$ contains a subword $e_i\rr{c}f_j,\ c\in\set{3,9}$
      or its reverse, it does not contain any subword
      $e_{i'}\rr{c'}f_{j'},\ c'\not\equiv0\pmod3,\ e_i\approx e_{i'}$
      (equivalently, $f_j\approx f_{j'}$) or its reverse. 
    \end{enumerate}
    Here the \emph{reverse} to the word $w$ is the word
    $w^*=x_l r_l x_{l-1}\dots x_2r_2x_1$. We call $l$ the
    \emph{length} of the word $w$.
    \item We define a \emph{cycle} as a pair $z=(w,r_1)$, where $w$ is a
    word such that $r_2=r_l=\sim$ and $r_k\ne\rr{c}$ with $v(c)=2$,
    while $r_1=\rr{c}$ with $v(c)\ne2$ and $x_lr_1x_1$ in $\dX$. For
    such a cycle we set $x_{ql+k}=x_k$ and $r_{ql+k}=r_k$ for any $q$
    and $1\le k\le l$.
     \item  The \emph{$m$-th shift} of the cycle $z=(w,r_1)$ is defined
     as the cycle $z^{(m)}=(w^{(m)},r_{2m+1})$, where
     $w^{(m)}=x_{2m+1}r_{2m+2}x_{2m+2}\dots r_{2m}x_{2m}$.
     \item  A cycle $(w,r_1)$ is called \emph{periodic} if $w$ is
     of the form $w=vr_1vr_1\dots r_1v$ for a shorter cycle $(v,r_1)$.  
    \item  We call two words, $w$ and
    $w'=x_1r'_2x_2r'_3\dots r'_lx_l$ (with the same $x_k$),
    \emph{elementary congruent} if there are two indices $k_1,k_2$ such
    that 
    \begin{align*}
      r_{k_1}&=\,\rr{3c},\ r_{k_2}=\,\rr{d}\ \text{ for some }\
      c\ne3,\,d\ne3,\\ 
      r'_{k_1}&=\,\rr{c},\  r'_{k_2}=\,\rr{3d},\\
      r'_k&=r_k\ \text{ for }\ k\notin\set{k_1,k_2},\\
      x_{k_1}&\approx x_{k_2}\ \text{ or }\ x_{k_1}\approx x_{k_2-1}.
    \end{align*}
     \item  We call two words $w,w'$ \emph{congruent}
     and write $w\equiv w'$ if there is a sequence of words
      $w=w_1,w_2,\dots,w_n=w$ such that $w_k$
     and $w_{k+1}$ are elementary congruent for $1\le k<n$. We call two
     cycles $z=(w,r_1)$ and $z'=(w',r'_1)$ \emph{congruent} and write
     $z\equiv z'$ if $w'\equiv z$ and $r_1'=r_1$.
    \end{enumerate}     
  \end{defin}

   We recall that two polyhedra $X,Y$ are called \emph{congruent} if
   $X\vee Z\simeq Y\vee Z$ for some polyhedron $Z$. Then we write
   $X\equiv Y$.

   \begin{theorem}\label{main}
     \begin{enumerate}
      \item  Every word $w$ defines an indecomposable polyhedron $P(w)$ from
       $\cT_7$, called \emph{string polyhedron}.
      \item  Let $\pi(t)\ne t$ be an irreducible polynomial over the
     field $\mZ/2$. Every triple $(z,\pi(t),m)$, where is a
     non-periodic cycle and $m\in\mN$, defines an indecomposable
     polyhedron $P(z,\pi,m)$ from $\cT_7$, called \emph{band
     polyhedron}.
      \item  Every indecomposable polyhedron from $\cT_7$ is congruent
      either to a string or to a band one. 
      \item  $P(w)\equiv P(w')$ \iff either $w'\equiv w$ or $w'\equiv w^*$.
      \item  $P(z,\pi(t),m)\equiv P(z',\pi'(t),m)$ \iff $m=m'$ and one
     of the following possibilities hold:
     \begin{enumerate}
      \item  $\pi'(t)=\pi(t)$ and either $z'\equiv z^{(k)}$ with $k$
      even or $z'={z^*}^{(k)}$ with $k$ odd;
      \item  $\pi'(t)=t^d\pi(1/t)$, where $d=\deg\pi$,  and either
      $z'=z^{(k)}$ with $k$ odd or $z'={z^*}^{(k)}$ with $k$ even.
     \end{enumerate}
      \item  Neither string polyhedron is congruent to a band one.
      \item  
     \end{enumerate}
   \end{theorem}

    The cofibration sequence 
 \[
   A\xarr f B \to Cf \to SA,\quad A\in\cA^0,\ B\in\cB^0,
 \]
    and the attachment diagram of a string polyhedron $P(w)$ is constructed as
    follows. 
    \begin{enumerate}
     \item  The indecomposable summands of $A$ correspond to the
     following subwords of $w$ (or their reverse):
     \begin{align*}
         S^{10}&\hspace*{-4em} &\text{to} &\hspace*{-4em}& f_1&\sim f'_1,\\
         S^{11}&\hspace*{-4em} &\text{to} &\hspace*{-4em}& f_2&\sim f'_2,\\
         S^{12}&\hspace*{-4em} &\text{to} &\hspace*{-4em}& f_3&\sim f'_3,\\
         C^{12}&\hspace*{-4em} &\text{to} &\hspace*{-4em}& f_4&\sim f_5.
     \end{align*}
      \item   The indecomposable summands of $B$ correspond to the
     following subwords of $w$ (or their reverse):
     \begin{align*}
         S^{7}&\hspace*{-4em} &\text{to} &\hspace*{-4em}& e_1,\\
         C_2^{10}&\hspace*{-4em} &\text{to} &\hspace*{-4em}& e_2&\sim e_3,\\
         C^9&\hspace*{-4em} &\text{to} &\hspace*{-4em}& e_4&\sim e_5,\\
         S^8&\hspace*{-4em} &\text{to} &\hspace*{-4em}& e_6&\sim
         e'_6,\\
     C^{10}&\hspace*{-4em} &\text{to} &\hspace*{-4em}& e_7,\\
     S^9&\hspace*{-4em} &\text{to} &\hspace*{-4em}& e_9\sim e'_9,\\
     S^{10}&\hspace*{-4em} &\text{to} &\hspace*{-4em}& e_{10}.
     \end{align*}
      \item  The attachments correspond to the subwords $e_i\rr{c}f_j$
      (or their reverse). Namely, such an attachment starts at the
      $f$-end of the corresponding subword and ends at its $e$-end;
      the number $c$ shows which multiple of the generator of the
      $(ij)$-th group from Table~9 must be taken. 
    \end{enumerate}

  For instance, if
  \begin{multline*}
    w=e_{10}\rr1 f'_2\sim f_2\rr8 e_6\sim e_6'\rr1 f_1'\sim f_1\rr2
    e_4 \\
    \sim e_5\rr6 f_5\sim f_4 \rr1 e_2\sim e_3\rr1 f'_3\sim f_3\rr2
    e_5\\
    \sim e_4\rr3 f_1\sim f'_1\rr1 e'_9\sim e_9\rr{12}f_3\sim f'_3,
  \end{multline*}
 the polyhedron $P(w)$ has the attachment diagram
 
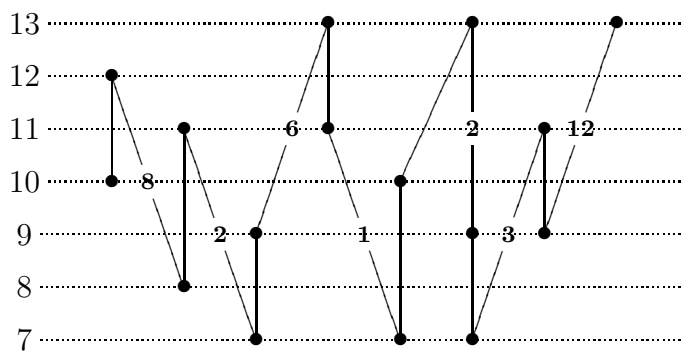
\[
   \xymatrix@R=1ex@C=1.7em{
   13 \ar@{.}[rrrrrrrrr] &&&& \cel \ar@{-}[dd]
   \ar@{-}[ddddl]|*+\txt{\scs 6} && \cel \ar@{-}[dddd]|*+\txt{\scs2} 
     \ar@{-}[dddl] && \cel \ar@{-}[ddddl]|*+\txt{\scs12}& *=0{}\\
   12 \ar@{.}[rrrrrrrrr] & \cel \ar@{-}[dd]
   \ar@{-}[ddddr]|*+\txt{\scs8} &&&&&&& & *=0{}\\ 
   11 \ar@{.}[rrrrrrrrr] && \cel \ar@{-}[ddd]
   \ar@{-}[ddddr]|*+\txt{\scs2} && \cel \ar@{-}[ddddr]|*+\txt{\scs1}
   &&& \cel 
   \ar@{-}[ddddl]|*+\txt{\scs3} \ar@{-}[dd] && *=0{}\\
   10 \ar@{.}[rrrrrrrrr] & \cel &&&& \cel \ar@{-}[ddd] &&&& *=0{}\\
   9 \ar@{.}[rrrrrrrrr]  &&&\cel \ar@{-}[dd] &&& \cel \ar@{-}[dd] & \cel && *=0{}\\
   8 \ar@{.}[rrrrrrrrr] && \cel &&&&&&& *=0{}\\
   7 \ar@{.}[rrrrrrrrr] &&& \cel && \cel & \cel  &&& *=0{}	} 
 \]
 
 Let now $P(z,\pi(t),m)$ be a band polyhedron. Replacing $w$ by $w^*$,
  we may suppose that $x_1\in\dE,\,x_n\in\dF$. Let also $\Phi$ be the
  Frobenius matrix with the characteristic polynomial
  $\pi(t)^m$. Then the cofibration sequence and the attachment diagram
  are constructed as follows.
  \begin{enumerate}
     \item  Do the construction as above for the word $w$.
     \item  Replace every summand $A_j$ of $A$ and every summand $B_i$
     of $B$ by $m$ copies of it, $A_{j1},\dots,A_{jm}$ and
     $B_{i1},\dots,B_{im}$.  
     \item  If there was an attachment $A_j\xarr c B_i$, replace it by
     the attachments $A_{jk}\xarr cB_{ik}\ (1\le k\le m)$.
     \item  If $A_j$ is the last summand of $A$, $B_i$ is the first 
     summand of $B$ and $r_1=\,\rr{c}$, add new attachments
     $A_{jk}\xarr cB_{il}$ in all cases, when the $(lk)$-th
     coefficient of the matrix $\Phi$ is nonzero.
  \end{enumerate}

  For instance, consider the band polyhedron $P(z,t^2+t+1,3)$
  $z=(w,\rr1)$, where
 \[
  w=e_2\sim e_3\rr1 f'_3\sim f_3\rr2 e_9\sim e'_9\rr1 f'_1\sim f_1\rr3
  e_4\sim e_5 \rr6 f_5\sim f_4.
 \]
  Then the attachment diagram is
 
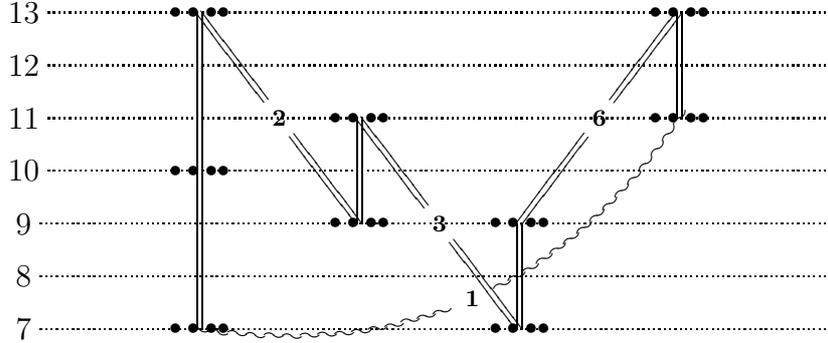
\[
  \xymatrix@R=1ex@C=4.5em{
 13 \ar@{.}[rrrrr] &\clls \ar@{=}[ddd] \ar@{=}[ddddr]|*+\txt{\scs2} & &
 &\clls \ar@{=}[dd] \ar@{=}[ddddl]|*+\txt{\scs6} &  *=0{}\\ 
 12 \ar@{.}[rrrrr] & & & & &  *=0{}\\
 11 \ar@{.}[rrrrr] & &\clls \ar@{=}[dd] \ar@{=}[ddddr]|*+\txt{\scs3} &
 &\clls &  *=0{}\\ 
 10 \ar@{.}[rrrrr] &\clls \ar@{=}[ddd] & & & &  *=0{}\\
  9 \ar@{.}[rrrrr] & &\clls &\clls \ar@{=}[dd] & &  *=0{}\\
  8 \ar@{.}[rrrrr] & & & & &  *=0{}\\
  7 \ar@{.}[rrrrr] &\clls \ar@{~}@/_6ex/[uuuurrr]|*+\txt{\scs1} & &\clls & &  *=0{}
  }
 \]
 Here the double lines show the attachments like
 
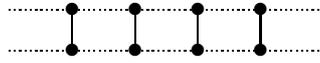
\[
  \xymatrix@=3ex@C=2em{ *=0{} \ar@{.}[rrrrr] & \cel \ar@{-}[d]  & \cel
  \ar@{-}[d]  & \cel \ar@{-}[d]  & \cel \ar@{-}[d] & *=0{}\\
 *=0{} \ar@{.}[rrrrr] & \cel   & \cel
 & \cel  & \cel  & *=0{}  }
 \]
 while the wavy line shows the attachment
 
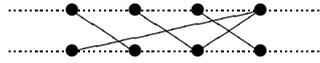
\[
  \xymatrix@=3ex@C=2em{ *=0{} \ar@{.}[rrrrr] & \cel \ar@{-}[dr]  & \cel
  \ar@{-}[dr]  & \cel \ar@{-}[dr]  & \cel \ar@{-}[dl] \ar@{-}[dlll] & *=0{}\\
 *=0{} \ar@{.}[rrrrr] & \cel   & \cel
 & \cel  & \cel  & *=0{}  }
 \]
 ruled by the Frobenius matrix with the characteristic polynomial
 $\pi(t)^2=t^4+t^2+1$, namely,
 \[\mtr{0&0&0&1\\1&0&0&0\\0&1&0&1\\0&0&1&0}.\]

  \section{Wild cases}
 \label{s7}

 Since we are dealing with additive categories that are not categories
 over a filed, we have to precise the notion of wildness. The
 following one seems to work in all known cases.

 \begin{defin}\label{wild}
   We call an additive category $\cC$ \emph{wild} if, there is a field
   $\Mk$ such that for every finitely generated $\Mk$-algebra $\La$
   there is a full subcategory $\cC_\La\sbe\cS$ and an epivalence
   $\cC_\La\to\La\Mod$ (the category of $\La$-modules that are finite
   dimensional over $\Mk$).
 \end{defin}

 One can see that for algebras over a field this definition is
 equivalent to the usual one (see, for instance \cite{d0}). One can
 also easily show that if a category $\cD$ is wild and there is an
 epivalence $\cC'\to\cD$ for a full subcategory $\cC'\sbe\cC$, then
 $\cC$ is wild as well.

 Now we present the results on wildness of categories $\cS_n$ and
 $\cT_n$. 

\begin{theorem}[Baues \cite{bd3}]\label{wild-1}
  If $n>4$, the category $\cS_n$ is wild.
\end{theorem}
\begin{proof}
  Obviously, one only has to prove the claim for $n=5$.
  The category $\cS_5$ contains the full subcategory $\cC=\cA\dag\cB$,
  where $\cA$ consists of bouquets of suspended Moore atom $A=M^6(2)$
  and $\cB$ consists of bouquets of suspended Moore atoms
  $B=M^8(2)$.Let $\cV={}_\cA\cS_\cB$. Since $\hos(B,A)=0$,
  Corollary~\ref{12} is applicable. Moreover, the ideal $\cJ$ in this
  case is zero, so $\cC/\cI\simeq\bim(\cV)$ with $\cI^2=0$, hence, the
  natural functor $\cC\to\bim(\cV)$ is an epivalence.

 Consider the cofibration sequence
 \begin{equation}\label{seq}
  S^7\xarr 2 S^7 \to A \to S^8 \xarr 2 S^8.
 \end{equation}
 Apply to it the functors $\hos(\_\,,S^6)$ and $\hos(\_\,,S^5)$. Taking
 into account the Hopf map $\eta:S^6\to S^5$ we get the commutative
 diagram with exact rows
 \[
 \begin{CD}
   0 @>>> \mZ/2 @>>> \hos(A,S^6) @>>> \mZ/2 @>>> 0 \\
   && @V\eta_*VV    @VVV  @VV\wr V \\
  0 @>>> \mZ/2 @>>> \hos(A,S^5) @>>> \mZ/2 @>>> 0\,.   
 \end{CD}
 \]
 Since $\eta^3=4\nu$, where $\nu$ is the element of order $8$ in
 $\hos(S^8,S^5)$, the map $\eta_*$ in this diagram is zero, therefore,
 the lower exact sequence splits and
 $\hos(A,S^5)\simeq\mZ/2\+\mZ/2$.Quite similarly, one shows that
 $\hos(S^8,B)\simeq\mZ/2\+\mZ/2$. Now apply the functors
 $\hos(\_\,,S^5)$ and $\hos(\_\,,B)$ to the exact sequence \eqref{seq}
 and take into account the map $S^5\to B$ form the definition of
 $B=M^6(2)$. Since $\hos(S^7,B)\simeq\mZ/2$, we get the commutative
 diagram with exact rows
 \[
 \begin{CD}
   0 @>>> \mZ/2 @>>> \hos(A,S^5) @>>> \mZ/2 @>>> 0 \\
   && @VVV    @VVV  @VV\wr V \\
  0 @>>> \mZ/2\+\mZ/2 @>>> \hos(A,B) @>>> \mZ/2 @>>> 0\,.   
 \end{CD}
 \]
 We know that the upper row of this diagram splits. Hence, the lower
 row splits too, so $\hos(A,B)\simeq(\mZ/2)^3$. Recall that
 $\es(A)\simeq\es(B)\simeq\mZ/4$ (Proposition~\ref{33}). Hence, there
 is an epivalence $\bim(\cV)\to\La\Mod$, where $\La$ is the path
 algebra of the quiver 
     $\xymatrix{\bullet \ar@/^/[r]\ar[r]\ar@/_/[r]&\bullet}$
 over the field $\mZ/2$. The
 latter is well-known to be wild, therefore, so is also $\cS_5$.
\end{proof}

\begin{theorem}[\cite{d2}]\label{wild-2}
  The category $\cT_n$ is wild for $n>7$.
\end{theorem}
\begin{proof}
  Again we only have to prove it for $n=8$. The category $\cT_8$
  contains the full subcategory $\cC=\cA\dag_{\cV}\cB$, where $\cA$ consists of
  bouquets of Chang atoms $C^14_2$ and $\cB$ consists of bouquets of
  spheres $\S^8$ and $S^11$, and
  $\cV={}_\cB\!\cS^0_{8,3}\!\cA$. Moreover,
  $\cI^0_{8,3}\cap\bim{\cV}=0$, so there is an epivalence
  $\cC\to\bim(\cV)$. Consider the cofibration sequence
  \[
    S^{13} \xarr{\eta^2} S^{11} \to C^{14}_2 \to S^{14} \to S^{12}
  \]
 and apply to it the functor $\hos(\_\,,S^{11})$. We get the exact
 sequence 
 \[
  \mZ/2 \xarr{(\eta^2)^*} \mZ/24 \to \hos(C_2^{14},S^{11}) \to \mZ \to
  \mZ/2, 
 \]
 wherefrom $\cS^0(C_2^{14},S^11)\simeq\mZ/12$. Moreover, there is a
 commutative diagram of cofibration sequences
 \[
 \begin{CD}
      S^{13} @>{\eta}>> S^{12} @>>> C^{14} @>>> S^{14} @>\eta>>S^{13}
      \\ 
      @V\id VV  @V\eta VV  @VVV  @VV\id V  @VV\eta V  \\
   S^{13} @>{\eta^2}>> S^{11} @>>> C^{14}_2 @>>> S^{14} @>>{\eta^2}> S^{12}\,. 
 \end{CD}
 \]
 Applying the functor $\hos(\_\,,S^8)$, we get the commutative diagram
 with exact rows
 \[
 \begin{CD}
  0 @>>> \mZ/2 @>>> \hos(C^{14}_2,S^8) @>>> \mZ/24 @>>> 0 \\
  && @VVV @VVV @VVV  \\
  0 @>>> \mZ/2 @>>> \hos(C^{14},S^8) @>>> 0 @>>> 0\,.
 \end{CD}
 \]
 (Recall that $\pis_{d+4}(S^d)=\pis_{d+5}(S^d)=0$ and
 $\pis_{d+6}(S^d)=\mZ/2$ \cite{to}). Therefore,
 $\cS^0(C_2^{14},S^8)\simeq\mZ/24\+\mZ/2$. So we present maps
 $a\in\cV(A,B)$, where $A\in\cA,\,B\in\cB$, as block-triangular
 matrices 
 \[
  a=\mtr{a_1&a_2\\0&a_3},
 \]
 where $a_1$ is with the coefficients from $\mZ/24$, $a_2$ is with
 coefficients from $\mZ/2$ and $a_3$ with coefficients from
 $\mZ/12$. On the other hand, maps $\al:A\to A'$, where $A,A'\in\cA$,
 and $\be:B\to B'$, where $B,B'\in\cB$ can be presented by block- triangular
 matrices
 \[
  \al=\mtr{\al_1&\al_2\\0&\al_3}\ \text{ and }\
  \be=\mtr{\be_1&\be_2\\0&\be_3},
 \]
 where $\al_2$ has coefficients from $\mZ/12$, $\be_2$ has
 coefficients from $\mZ/24$, other blocks have components from $\mZ$
 and $\al_1\equiv\al_3\mod2$. 

 We consider the full subcategory $\cC\sb\bim(\cV)$ consisting of all
 maps $a$ such that the corresponding blocks $a_1,a_2,a_3$ are of the
 form
 \[
  a_1=\mtr{6I&0&0\\0&12&0},\quad a_2=\mtr{0&I&0\\0&0&u} \quad
  a_3=\mtr{6v_1&6v_2&0}, 
 \]
 where the entries $I$ stand for identity matrices (not necessary of
 the same dimensions) and $u,v_1,v_2$ are arbitrary matrices with
 coefficients from $\mZ/2$ of proper sizes. We write $a=a(u,v_1,v_2)$. One can verify
 that if $(\al,\be)$ is a morphism $a(u,v_1,v_2)\to a(u',v_1',v_2')$,
 there are integral matrices $\ga_1,\ga_2,\ga_3$ such that
 $v_i\ga_1=\ga_2v_i\ (i=1,2)$ and $u\ga_3=\ga_1u$. Conversely, any
 given triple $\ga_1,\ga_2,\ga_3$ with these properties can be
 accomplished to a morphism $a(u,v_1,v_2)\to a(u',v_1',v_2')$. It
 gives rise to an epivalence $\cC\to\La\Mod$, where $\La$ is the path
 algebra of the quiver
 $\xymatrix{\bullet\ar[r]&\bullet\ar@/^/[r]\ar@/_/[r]&\bullet}$. It is
 known to be wild. Therefore, $\cT_8$ is wild as well.
\end{proof}

 \end{document}